%% file: paper.tex
\documentclass[final,3p,times]{elsarticle}

\usepackage[normalem]{ulem}

\usepackage{amssymb}
\usepackage{amsmath}
\usepackage{amsthm}
\usepackage{subeqnarray}
\usepackage{mathrsfs}
\usepackage{rotating}
\usepackage{enumitem}
\setlist[enumerate]{itemsep=1pt, topsep=2pt}
\setlist[itemize]{itemsep=1pt, topsep=2pt}
\usepackage{bm}             
\usepackage{setspace}
\usepackage{numprint}

\usepackage{algorithmic}
\usepackage[linesnumbered,ruled]{algorithm2e}

\usepackage[breaklinks,bookmarks=false]{hyperref}
\hypersetup{colorlinks, linkcolor=blue, citecolor=blue,
  urlcolor=blue}

\usepackage{cleveref}
\usepackage{amstext}

\biboptions{sort&compress}

\input{def.tex}

\usepackage[usenames]{color}

\newcommand{\Smaz}[1]{{}}

\usepackage[scientific-notation=true]{siunitx} 
\newcommand{\numf}[1]{\num[round-mode=places, scientific-notation=false, round-precision=0]{#1}}

\begin{document}


\begin{frontmatter}



  \title{Adaptive domain decomposition method for time-dependent problems with applications in 
    fluid dynamics} 

\author[Prague]{V{\'\i}t Dolej{\v s}{\'\i}}
\ead{vit.dolejsi@matfyz.cuni.cz}

\author[MUAV]{Jakub {\v S}{\'\i}stek}
\ead{sistek@math.cas.cz}

\address[Prague]{Charles University, Faculty of Mathematics and Physics,
Sokolovsk\'a 83, 186 75 Praha, Czech Republic}

\address[MUAV]{Institute of Mathematics of the Czech Academy of Sciences, 
Žitná 25, 110 00 Praha, Czech Republic}


\begin{abstract}
  We deal with the numerical solution of the time-dependent partial differential equations using the adaptive space-time discontinuous Galerkin (DG) method. The discretization leads to a nonlinear algebraic system at each time level, the size of the system is varying due to mesh adaptation. A Newton-like iterative solver leads to a sequence of linear algebraic systems which are solved by GMRES solver with a domain decomposition preconditioner. Particularly, we consider additive and hybrid two-level Schwarz preconditioners which are efficient and easy to implement for DG discretization. We study the convergence of the linear solver in dependence on the number of subdomains and the number of element of the coarse grid. We propose a simplified cost model measuring the computational costs in terms of floating-point operations, the speed of computation, and the wall-clock time for communications among computer cores. Moreover, the cost model serves as a base of the presented adaptive domain decomposition method which chooses the number of subdomains and the number of element of the coarse grid in order to minimize the computational costs. The efficiency of the proposed technique is demonstrated by two benchmark problems of compressible flow simulations.
\end{abstract}
\end{frontmatter}
  


\section{Introduction}

Numerical solution of time-dependent partial differential equations is
an important tool in many areas of science and engineering. Among various numerical techniques,
the {\em space-time discontinuous Galerkin} (STDG) method presents a very powerful discretization
approach. It is based on piece-wise polynomial
discontinuous approximation with respect to spatial and temporal variables.
The main advantages of STDG method are stability, accuracy and robustness,
the ability to deal with $hp$-mesh adaptation,
the use of arbitrary non-nested meshes for time-dependent problems,
and a relative simplicity of constructing domain decomposition preconditioners. We refer to
earlier papers \cite{KVV:NS1,KVV:NS2,vegt_sudirham_damme,RhebergenALL_JCP13} dealing with numerical
solution of fluid dynamics problems,
\cite{SmearsSuli_NM16} using STDG method for the solution of Hamilton--Jacobi--Hellman equations,
\cite{CangianiALL_SISC17} analysing STDG method for parabolic problems on prismatic meshes,
\cite{Wieners_CAMWA23} for solving linear transport problem,
or \cite{CoralloALL_JSC23} dealing with symmetric Friedrichs systems.
Papers \cite{Munz07,Dumbser2011,TAVELLI2018386,RomeoDumbserTavelli_CAMS21} also present related
research.
Main theoretical results can be found in \cite[Chapter~6]{DGM-book}.

Numerical solution of time-dependent problems typically involves re-meshing during
the evolution process in order to balance the accuracy and efficiency of the simulation.
In \cite{AMAtdp,Dolejsi_IJCFD24}, we developed
the {\em anisotropic $hp$-mesh adaptation} technique which allows to adapt 
the  size, shape, and orientation of triangular elements 
together with varying polynomial approximation degrees.
This adaptive technique minimizes the number of degrees of freedom so that
an interpolation error estimate is under the given tolerance.
The STDG method naturally handles grids generated by
the anisotropic mesh adaptation method since it is
a one-step method, and the approximate solutions at different time levels
are connected only weakly by a time penalty term.
The computational benefit of anisotropic mesh adaptation
for unsteady flow problems was presented in many papers, e.g., 
\cite{BelmeDervieuxAlauzet_JCP12,AlauzetLoseille_JCP18,FreyAlauzet,HabashiAll,SauvageALL_JCP24} and
the references cited therein.

The significant disadvantage of STDG method is the size of the arising nonlinear algebraic systems.
Particularly, using the approximation of degree $q\ge 0$, STDG method gives a system of
size $(q+1)\times \DoF$ at each time step where $\DoF$ is the number of spatial degrees of freedom.
The temporal accuracy of this scheme at the nodes of the time partition is $2q+1$.
On the other hand, the popular diagonally implicit Runge-Kutta (DIRK) methods
(\cite{DIRK,DIRK2,SODE2,DIRK3,DIRK4,LevyMay_CF25}) having
order $r\ge1$
require the solution of $r$ systems of size $\DoF$ at each time step which is typically
faster.
Therefore, to make the STDG approach competitive, efficient algebraic solvers are required.

The iterative solution of nonlinear algebraic systems leads to sequences of linear systems
that are frequently solved iteratively by Krylov methods \cite{Liesen20131} with suitable
preconditioners. A prominent role is played by the {\em domain decomposition} (DD) techniques
\cite{ToselliWidlund-DD05,Nataf-DDM} since they can employ the power of parallel computers.
Their combination with discontinuous Galerkin discretization leads to
a very comfortable setting.
Namely, non-overlapping domain decomposition preconditioners
do not require any special choice of the interface conditions among subdomains.
Furthermore, two-level Schwarz preconditioners, which accelerate the transfer
of information through the whole system,  can be constructed in a natural way, see e.g.,
\cite{Karashian2001,AntoniettiHouston_JSC11,AntoniettiALL_IJNAM16,Karakashian_IMA17}.

In this paper, we employ a domain decomposition preconditioners for the
solution of systems arising from adaptive STDG discretization of time-dependent partial differential
equations describing the motion of fluids. 
In particular, we adopt additive and hybrid two-level Schwarz techniques with two domain
decomposition parameters: the number of subdomains $M$ and the number of elements of the coarse grid
$n_0$. We note that the subdomains arise from re-grouping of the  mesh elements into
$M$ connected subdomains by standard tools like METIS \cite{metis} such that the number
of degrees of freedom of each subdomain is equilibrated. Moreover, the coarse mesh arises from
the splitting of each subdomains into $s$ coarse macro-elements, hence the coarse mesh
has $n_0=s\times M$ elements.

It is known and documented also in this paper, that larger $M$ allows to use more parallelization
but decreases the speed of convergence. Similarly, a finer coarse mesh provides more
information and accelerates the convergence but the corresponding system is larger and, therefore,
more expensive to solve. In this paper we consider the following question: How to choose
the number of subdomains $M$ and the coarse mesh factor $s$ in an optimal way?

In order to easily analyse the computational costs, we introduce a new
{\em computational cost model}
measuring the number of floating point operations ({\flops}) per one core,
the performance (speed of computation), and the wall-clock time of communications
among computer cores.
We study numerically, how these costs depend on the domain decomposition parameters
$M$ and $s$. This is the first novelty of this paper.
Moreover, based on the
computational cost model, we propose an adaptive technique which estimate the computational costs
a priori and sets
the optimal parameters $M$ and $s$ for domain decomposition. 
The efficiency of this technique is demonstrated by several numerical experiments.
The proposed framework can be applied to an arbitrary adaptive method solving any time-dependent
problem. We emphasise that the goal of this paper is not to study and demonstrate the
scalability of the preconditioners, but to develop a strategy for choosing optimally the
domain decomposition parameters for a particular problem.

In Section~\ref{sec:problem}, we present the model describing the motion of viscous
compressible fluid, its discretization by the STDG method with several solution strategies is
briefly given in Section~\ref{sec:DGM}. The domain decomposition (DD) based preconditioners
are explained in Section~\ref{sec:DDM} where we also formulate the fundamental
problem of this paper. In Section~\ref{sec:costs}, we discuss the computational costs
of the DD preconditioners. The dependence of the computational costs on
DD parameters $M$ and $s$ is studied numerically.
Moreover, we present the resulting adaptive domain decomposition technique in
Section~\ref{sec:ADDM}, its efficiency is demonstrated by two examples
in Section~\ref{sec:numer}, and we conclude with several remarks in Section~\ref{sec:concl}.


\section{Governing equations}
\label{sec:problem}

We consider the system of the compressible Navier-Stokes equations with state equation
for perfect gas.
Let $\Om\subset\R^2$ be the computational domain occupied by the fluid,
$T>0$ the physical time  to be reached,
and $Q_T:= \Om\times(0,T)$.
Moreover, the symbols $\pdt{}$ and $\pdd{}{x_i}$ denote the partial derivatives with respect to
$t$ and $x_i,\ i=1,2$, respectively.
Then the motion of the fluid is described by the system of $n=4$ convection-diffusion equations
\begin{align} \label{eq:NS}
  \pdt{\w} + \sum\nolimits_{i=1}^2 \pdd{}{x_i} \left(\bF_i(\w)  -  \bR_i(\w, \nabla \w)  \right)
  = \bg(\w) \quad\mbox{in} \ Q_T,
\end{align}
where 
$\w=\w(x,t): Q_T\to \R^n$ is the unknown state vector,
$\bF_i: \R^n \to \R^n$, $i=1,2$ represent the convective fluxes,
$\bR_i: \R^n \to \R^n$, $i=1,2$ represent the diffusive fluxes,
and $\bg: Q_T\to\R^n $ denotes the exterior forces.
System \eqref{eq:NS} is accompanied by the initial and boundary conditions.
Particularly, we set $\w(x,0) = \w_0(x)$ in $\Om$,
where $\w_0$ is a given function.
Moreover, 
various problem-dependent boundary conditions are prescribed on the boundary $\gom:=\pd\Om$.
The state vector and the fluxes read
\begin{align}
  \label{eq:w}
  \w &=(\rho,\,\rho v_1,\rho v_{2},\,\ener)^{\TT},\qquad
  \bF_i = (\rho v_i,\,\rho v_i v_1 + \delta_{i1}\press,\,\rho v_i v_2 + \delta_{i2}\press,\,
  (\ener+\press)\,v_i)^{\TT},\qquad i=1,2, \\
   \bR_i &= \bigg(0,\,\stress_{i1},\,\stress_{i2},
   \sum\nolimits_{k=1}^2 \stress_{ik}v_i
   + \frac{\tilde{\mu}\,\ccP}{\Pr }
   \pdd{\temp}{x_i}\bigg)^{\TT},\ i=1,2, \qquad
  \bg = \left(0,\, -\rho g k_1,\, -\rho g k_2,\, -\rho g \bk\cdot\bkv\right)^{\TT}, \notag
\end{align}
where
$\rho$ is the density,
$\press$ is the pressure, 
$\vv=(v_1,v_2)^\TT$ is the velocity vector, 
$\temp$ is the absolute temperature and
$\ener = \rho\ccV\temp + \rho |\bkv|^2/2$ is the energy per unit volume
including the interior and kinetic energies (and excluding the gravitational energy).
Moreover, 
$\ccV >0$ and 
$\ccP > 0$ are the specific heat capacities at constant volume and pressure, respectively, 
$\kappa= \ccP/\ccV > 1$ is the  Poisson adiabatic constant,
$R=\ccP-\ccV $ is the gas constant,
$\tilde{\mu}$
is the dynamic viscosity, $g=9.81\,\mathrm{m}\cdot\mathrm{s}^{-2}$ is the gravity constant,
$\bk=(k_1,k_2)^\TT$ is the upward pointing unit vector, $\Pr$ is the Prandtl number,
and $\stress_{i j},\ i,j=1,2$ denote the components of the viscous part of the stress tensor.
Symbol $\delta_{i j}$ is the Kronecker delta and the symbol $\cdot$ denotes the scalar product.
The relations \eqref{eq:NS} -- \eqref{eq:w} are accompanied by the 
constitutive relations of the perfect gas and the viscous part of the stress tensor as
\begin{align}
  \label{eq:tau}
  \press=R \rho \theta \qquad\mbox{and} \qquad
  \stress_{ij} = \tilde{\mu}  \left( 
 \pdd{v_i}{x_j} + \pdd{v_j}{x_i}\right)
 - 
 \tfrac23  \delta_{ij} 
 \left(\pdd{v_1}{x_1} + \pdd{v_2}{x_2}\right),
 \quad i,j=1,2, \quad \mbox{respectively.}
\end{align}

The efficiency of the presented domain decomposition techniques is demonstrated also on
a two-dimensional non-hydrostatic mesoscale atmospheric modeling problem where
the preferred physical quantity is  the potential temperature $\tempP$ given as
\begin{align}
  \label{eq:tempP}
  \tempP = {\temp}/{\pressE},\qquad\mbox{where}\quad 
  \pressE = \left({\press}/{\press_0}\right)^{(\kappa-1)/\kappa}
\end{align}
is the Exner pressure, and $\press_0 = 10^5\,\si{Pa}$ is the reference pressure.
The initial condition is typically set as a perturbation in terms of the potential
temperature of a steady-state flow where the mean values $\bar{\tempP}$ and $\bar{\pressE}$
of the Exner pressure and the potential temperature 
are in a hydrostatic balance (cf.~\cite{GiraldoRostelli_JCP08}), i.e.,
\begin{align}
  \label{balance}
  \ccP \bar{\tempP} \frac{\dd \bar{\pressE}}{\dd x_2} = - g.
\end{align}


\section{Discontinuous Galerkin discretization}
\label{sec:DGM}

We discretize problem \eqref{eq:NS}  by the 
{\em space-time discontinuous Galerkin method}.
A detailed description of the method  can be found, e.g., in \cite[Chapters~8--9]{DGM-book}
or \cite{Dolejsi_IJCFD24}.
Here, we present only the final formulas for completeness.

\subsection{Function spaces}
\label{sec:dgm1}
We introduce a partition of $(0,T)$ by  $0=t_0<t_1<\ldots <t_r=T$ and set
$I_m =(t_{m-1},t_{m})$, and  $\tau_m = t_m- t_{m-1}$ for  $m=1,\ldots,r$.
For each $t_m,\ m=0,\ldots,r$, we consider a 
mesh $\Thm$ consisting of a finite number 
of closed 
triangles $K$ with mutually disjoint interiors and covering $\oO$.
The meshes $\Thm$ can differ for $m=0,\dots,r$.
For each $K\in\Thm$, $m=0,\dots,r$, we assign a polynomial approximation degree
with respect to space $p_K\ge 1$.
Moreover, let $q\ge0$ be the polynomial approximation degree with respect to time which is
constant for all elements.
We define the spaces of {\em discontinuous piecewise polynomial}
functions on the space-time layer $\Om\times I_m$ by
\begin{align}
  \label{bSm}
  \bSm & := \left\{ \bpsi:\Om\times I_m \to \R^n;\ \ 
  \bpsi(x,t)|_{K\times I_m} = \sum\nolimits_{j=0}^q t^j \bvas(x),\
  \bvas|_K\in [P_{\pK}(K)]^n,\  K\in\Thm\right\},
  \qquad \ m=0,\dots,r, 
\end{align}
where $[P_{\pK}(K)]^n$ denotes the vector-valued space of all polynomials on $K$ of
degree $\leq \pK$, $K\in\Thm$. 

To proceed to an algebraic representation, a suitable basis of $\bSm$
has to be chosen.
In the context of discontinuous Galerkin discretization, it is natural to choose the
basis functions having support restricted to only one space-time element
$K\times I_m$, $K\in\Thm$, $m=1,\dots,r$.
Therefore, we consider the basis of the space $\bSm$,
$m=1,\dots,r$, 
\begin{align}
  \label{basis1}
  \Bm&:=\left\{\bpsi_{\balpha},\ 
  \balpha \in \BBm\right\},
\ \mbox{where}\ 
  \BBm:=\Big\{{\balpha};\ {\balpha}=\{K,k,i,j\},\ j=0,\dots,q,\
  i=1,\dots,\ddK,\ k=1,\dots,n,\ K\in\Thm \Big\}
\end{align}
is the set of multi-indices ${\balpha}$, index $i$ corresponds to
the spatial degrees of freedom $1,\dots,\ddK:=(\pK+1)(\pK+2)/2$ (= dimension of $P_{\pK}(K)$),
index $j$ corresponds to temporal degrees of freedom $0,\dots,q$, index $k=1,\dots,n$
denotes the equation
and $K\in\Thm$ is the mesh element.
Obviously, the cardinality of set $\BBm$ is equal to $\Nhm=\dim \bSm$, cf.~\eqref{Nhm}.
Moreover, we denote by $K_{\balpha}\in\Th$ the element defining multi-index ${\balpha}\in\BBm$.
Then any basis function
$\bpsi_{\balpha}$, ${\balpha}\in\Bm$ has support restricted only to element $K_{\balpha}\in \Thm$,
$m=1,\dots,r$.




\subsection{Space semi-discretization}
\label{sec:semi}

To derive the space discontinuous Galerkin (DG) semi-discretization,
we multiply  \eqref{eq:NS}, by $\bpsi\in\bSm$, integrate over
$K\in\Thm$, apply Green's theorem, sum over all $K\in\Thm$,
approximate the convective and diffusive fluxes through element boundaries,
and add stabilization terms vanishing for any smooth solution $\w$. Then, we obtain
an abstract identity
\begin{align}
  \label{ahm}
  \Lsp{\pdt{\w(\cdot,t)}}{\bpsi(\cdot, t)}
  +  \ahm\left(\w(\cdot,t),\bpsi(\cdot,t) \right) = 0
  \qquad \forall \bpsi\in\bSm,\quad t\in I_m,\ m=1,\dots,r,
\end{align}
where $\Lsp{\cdot}{\cdot}$ denotes the $L^2(\Om)$-scalar product and 
$\ahm$ 
is a form arising from DG discretization of \eqref{eq:NS}.
This form is linear with respect to its second argument.
We refer to \cite[Chapters 8--9]{DGM-book} or \cite{Dolejsi_IJCFD24} for the particular form
of $\ahm$.
We only note here that it 
contains integral over elements $K\in\Thm$ and their edges.

To treat the nonlinearity of system \eqref{ahm}, we introduce its linearized form
\begin{align}
  \label{ahL}
  \ahmL(\bbw, \w, \bpsi)
  \approx \frac{D}{D\w} \ahm(\bbw,\bpsi)
\end{align}
where $ \frac{D}{D\w} \ahm(\bbw,\bpsi)$ denotes the Fr\'echet derivative of
$\ahm$ with respect to $\w$ evaluated at $\bbw$.
It can be derived directly by the  differentiation of $\ahm$, or
a suitable approximation can be used too, see the references mentioned above.
Obviously, $\ahL$ is linear with respect to its second and third arguments.
%

\subsection{Full space-time discretization}
\label{sec:STDG}
To define the full space-time DG discretization of \eqref{eq:NS},
we integrate \eqref{ahm} over $I_m$ and add a ``time-stabilization'' terms vanishing
for solution continuous with respect to time.
Particularly, 
we introduce the
forms
\begin{align}
  \label{Ah}
   \Ahm(\w, \bpsi) & :=  \int_{I_m}\Bigl( \Lsp{\pdt{\w}}{\bpsi} + \ahm(\w,\psi)\Bigr)\dt
  + \LSP{\tjump{\w}_{m-1}}{\bpsi|_{m-1}^{+}}{\Om}, \\
  \AhmL(\bbw, \w, \bpsi) & :=  \int_{I_m}\left( \Lsp{\pdt{\w}}{\bpsi}
  + \ahmL(\bbw,\w,\psi)\right)\dt
   + \LSP{\w|_{m-1}^{+}}{\bpsi|_{m-1}^{+}}{\Om}, \quad \bbw,\w,\bpsi\in\bS,\quad
   m=1,\dots,r,\notag
\end{align}
where $\bS$ is the function space over the whole space-time cylinder $\Om\times(0,T)$ given by
\begin{align}
  \label{bS}
  \bS:=\{ \wht\in L^2(\Om\times (0,T)),\ \wht|_{\Om\times I_m} =: \whtm \in \bSm, m=1,\dots,r\}.
\end{align}
The symbol  $\tjump{\cdot}_{m-1}$ in \eqref{Ah} denotes the jump  with respect to time
on the time level $t_m,\ m=0,\dots, r-1$, namely
\begin{align}
  \label{tjump}
  \tjump{\bpsi}_m := \bpsi|_m^+ - \bpsi|_m^-,
  \qquad \bpsi|_m^\pm := \lim_{\epsilon\rightarrow 0\pm}\bpsi(t_m+\epsilon),\qquad \bpsi\in\bS,
\end{align}
where $\bpsi|_0^-$ is typically taken from the initial condition.
The term $\LSP{\tjump{\w}_{m-1}}{\bpsi|_m^{+}}{\Om}$ in \eqref{Ah} joins
together the solution on time intervals $I_{m-1}$ and $I_m$, and it 
represents the ``time penalty''.
This term is applicable also for $\bSmM \not=\bSm$ which is advantageous
for varying meshes during the adaptation process.
Form $\AhmL$ 
is linear with respect to its second and third arguments.
Due to 
the choice of basis \eqref{basis1}, term
$\AhmL(\bbw,\bpsi_{{\bbeta}},\bpsi_{{\balpha}})$ for $\bpsi_{{\balpha}},\bpsi_{{\bbeta}}\in\bSm$
is non-vanishing if $K_{\balpha} = K_{\bbeta}$ or $K_{\balpha}$ and $K_{\bbeta}$
have a common edge.
For a detailed derivation of time DG discretization, including the numerical analysis,
we refer to \cite[Chapter~6]{DGM-book}.

Now, we formulate the final numerical scheme.
\begin{definition}
  We say that the function $\wht=\{\whtm,\ m=1,\dots,r\}\in\bS$ (cf.~\eqref{bS})
  is the {\em space-time discontinuous Galerkin} (STDG) solution
  of \eqref{eq:NS} if
  \begin{align}
    \label{STDGM}
     \Ahm(\wht,\bpsi_h) = 0 \qquad \forall \bpsi_h \in\bSm, \quad m=1,\dots,r,
  \end{align}
  where we set $\wht|_0^- := \w_0$ (= the initial condition).
\end{definition}

\subsection{Solution strategy}

For any $m=1,\dots,r$, formula \eqref{STDGM} represents the system of $\Nhm$
algebraic equations, 
\begin{align}
  \label{Nhm}
  \Nhm = \dim\bSm =
  (q+1)\, n\, \sum_{K\in\Thm}\nolimits \ddK,  \qquad \ddK = (\pK+1)(\pK+2)/2,\ K\in\Thm,
\end{align}
for the solution $\whtm\in\bSm$ (depending also on $\whtmM$). Let $\Wm\in\R^{\Nhm}$ be
the algebraic representation of $\whtm$ in basis $\Bm$, cf.~\eqref{basis1}.
Then the algebraic representation of \eqref{STDGM} is written as 
\begin{align}
  \label{newton0}
  \Fm(\Wm) = 0,\qquad \mbox{where }
  \Fm:\R^{\Nhm} \to \R^{\Nhm},\ \Fm =\{ (\Fm)_{{\balpha}}\}_{{\balpha}\in\BBm},\quad
  (\Fm)_{{\balpha}} = \Ahm(\whtm, \bpsi_{{\balpha}}),
\end{align}
where $\bpsi_{{\balpha}}\in\Bm$ is the basis function.
Similarly, let matrix $\Am(\bW)$ be the algebraic representation of $\AhmL(\bbw,\cdot,\cdot)$,
where $\bW$ is the representation of $\bbw$. Namely
\begin{align}
  \label{newton1}
  \Am:\R^{\Nhm} \to \R^{\Nhm\times\Nhm},\quad
  \Am(\bW) = \left\{ (\Am(\bW))_{{\balpha},{\bbeta}} \right\}_{{\balpha},{\bbeta}\in\BBm},
  \quad
  (\Am(\bW))_{{\balpha},{\bbeta}} = \AhmL(\bbw, \bpsi_{{\bbeta}},\bpsi_{{\balpha}}).
\end{align}
In virtue of \eqref{ahL}, matrix $\Am$ is an approximation of the Jacobian of $\Fm$. 
Then the damped Newton-like method, computing an approximation of $\Wm$ by the sequence
$\Wml$, $\ell=0,1,\dots$, reads
\begin{subequations}
  \label{newton2}
  \begin{align}
    \label{newton2a}
    &\Wml := \WmlM - \lambda_\ell \dl,\qquad \ell=1,2,\dots,\\
    \label{newton2b}
    \mbox{where } \dl\in \R^{\Nhm} \mbox{ solves }   \quad
    &\Am(\WmlM) \dl = \Fm(\WmlM),
  \end{align}
\end{subequations}
the initial guess $\WmN$ is taken from the previous time layer, i.e., $\WmN:=\WmM$, and
$0< \lambda_\ell \leq 1$ is the damping parameter. It is chosen based on the monitoring
function $\zeta_{\ell}$ given by
\begin{align}
  \label{newton3}
  \zeta_\ell:= {\norm{\Fm(\Wml)}{2}}\,/\, {\norm{\Fm(\WmlM)}{2}},\qquad \ell=1,2,\dots,
\end{align}
where $\norm{\cdot}{2}$ is the Euclidean norm in $\R^{\Nhm}$.
For each iteration step $\ell$, we solve system \eqref{newton2b}
and perform step  \eqref{newton2a} with $\lambda_\ell=1$.
If $\zeta_\ell < 1$, this step is accepted and we proceed to $\ell+1$.
If $\zeta_\ell \ge 1$, we reduce $\lambda_\ell$ by a factor smaller than 1 (0.65 in our case)
and repeat step  \eqref{newton2a} until $\zeta_\ell < 1$.
The whole iterative process \eqref{newton2} is stopped when a suitable stopping criterion is
achieved, cf. Section~\ref{sec:solvers}.
In the following, we call \eqref{newton2} the {\em Newton method} for simplicity
although its various modification can be considered.
\begin{remark}
  \label{rem:refresh}
  The most time consuming parts of the computational process are the evaluation of entries
  of matrix $\Am$ in \eqref{newton2b} and, especially, the solution  of \eqref{newton2b}.
  To reduce the computational costs, we do not update matrix $\Am(\WmlM)$ for
  each $\ell$ since it serves only as an approximation of the Jacobian in the Newton method.
  In particular, if the monitoring function $\zeta_\ell < 1/2$, we keep the actual evaluation
  of $\Am(\cdot)$. Otherwise, matrix $\Am$ is refreshed in the actual approximation $\WmlM$.
  In our experience, this strategy slightly increases the number of iterations $\ell$
  to satisfy the stopping
  criterion but it saves the computational time since evaluation of the right-hand side
  $\Fm(\WmlM)$ is much cheaper than the evaluation of $\Am(\WmlM)$,
  and the factorization of the preconditioner (cf. Section~\ref{sec:DDM}) is also kept.
\end{remark}

At each step $\ell$ of the nonlinear iterative solver \eqref{newton2b}, we have to solve
the linear algebraic system which is discusses in Section~\ref{sec:DDM}.

\subsection{Stopping criteria for iterative solver and choice of the time steps}
\label{sec:solvers}

Finally, we briefly mention the remaining parts of the computational process, particularly, the
choice of the time steps $\tau_m=|I_m|$, $m=1,\dots,r$ in \eqref{Ah} and \eqref{STDGM}, and
the stopping criteria for nonlinear iterative process  \eqref{newton2}.
Obviously, the errors arising from
the space and time discretization have to be balanced, and the algebraic errors should
be below the discretization one.

Therefore, in virtue of \cite{Dolejsi_IJCFD24,stdgm_est}, for each time level $m=1,\dots,r$,
we define the estimators
\begin{align}
  \label{etas}
  \etaAm(\wht) &: = \sup_{0\not=\bpsi_h\in \bSm}
  \frac{\Ahm(\wht,\bpsi_h)}{\norm{\bpsi_h}{X}}, 
  \ \
  \etaSm(\wht) = \sup_{0\not=\bpsi_h\in \bSmhh}
  \frac{\Ahm(\wht,\bpsi_h)}{\norm{\bpsi_h}{X}},
  \ \ 
  \etaTm(\wht) = \sup_{0\not=\bpsi_h\in \bSmtt}
  \frac{\Ahm(\wht,\bpsi_h)}{\norm{\bpsi_h}{X}}, 
\end{align}
where $\bSm$ is given by \eqref{bSm} and
\begin{align}
  \label{etas2}
  \bSmhh & := \left\{ \bpsi:\Om\times I_m \to \R^n;\ \ 
  \bpsi(x,t)|_{K\times I_m} = \sum\nolimits_{j=0}^q t^j \bvas(x),\
  \bvas\in [P_{\pK+1}(K)]^n,\  K\in\Thm\right\},  \\
  \bSmtt & := \left\{ \bpsi:\Om\times I_m \to \R^n;\ \ 
  \bpsi(x,t)|_{K\times I_m} = \sum\nolimits_{j=0}^{q+1} t^j \bvas(x),\
  \bvas\in [P_{\pK}(K)]^n,\  K\in\Thm\right\}, \notag 
\end{align}
are its enrichment with respect to space and time. 
In addition,
$\norm{\bpsi}{X}=\left(
\normP{\bpsi}{L^2(\Om\times I_m)}{2} +
\normP{\nabla \bpsi}{L^2(\Om\times I_m)}{2} +
\normP{\pdt\bpsi}{L^2(\Om\times I_m)}{2}\right)^{1/2}$.\\[2pt]
We showed in \cite{stdgm_est} that
$\etaSm$ and $\etaTm$ correspond to the spatial and temporal discretization errors and
$\etaAm$ corresponds to the algebraic error.
Obviously, $\etaAm(\wht)=0$ for the exactly computed solution $\wht$ of \eqref{newton0}.
Hence, we stop the nonlinear iterative process \eqref{newton2} if the following condition
is valid
\begin{align}
  \label{CA}
  \etaAm(\wht) \leq \CA \min(\etaSm(\wht), \etaTm(\wht)),\qquad \mbox{ where } \CA\in(0,1).
\end{align}
Similarly, the time step $\tau_m$ is chosen such that
\begin{align}
  \label{CT}
  \etaTm(\wht) \approx \CT \etaSm(\wht), \qquad \mbox{ where } \CT \in(0,1).
\end{align}
For more details, we refer to \cite{Dolejsi_IJCFD24,stdgm_est}.

\subsection{Anisotropic mesh adaptation}
\label{sec:AMA}

We briefly mention the main idea of the {\em anisotropic $hp$-mesh adaptation method} ({\hpAMA})
developed in \cite{AMAtdp,Dolejsi_IJCFD24}. The idea is to define
a sequence of meshes $\Thm$ and spaces $\bSm$ for all $m=1,\dots, r$ such that
the interpolation error is under the given tolerance, and $\Nhm =\dim\bSm$ is minimal.
The interpolation error is estimated from a higher-order local reconstruction of the
quantity of interest (typically the density) by a least square method where
the available information from neighboring elements is employed.
If the  interpolation error estimate is below the tolerance, we use the given mesh
for the next time step, otherwise a re-meshing is carried out and the time step
has to be repeated.

The re-meshing techniques construct the new mesh and piecewise polynomial space by
minimizing $\Nhm$.  The anisotropic $hp$-mesh adaptation technique admits to
modify size, shape and orientation of mesh elements, as well as local polynomial
approximation degree with respect to space. The resulting meshes are non-matching,
non-nested but the solution on consequent time intervals are joined together by the time penalty,
cf. Section~\ref{sec:STDG}. For details, we refer to \cite{AMA-book}.

\section{Domain decomposition preconditioners}
\label{sec:DDM}

Now, we focus on the solution of the sequence of
linear algebraic systems \eqref{newton2b} which are solved iteratively by
the GMRES method with domain decomposition based preconditioners.
Therefore, instead of $\Am\bm{d}=\Fm$, we solve the preconditioned  system
$\NN_m^{-1}\Am\bm{d}=\NN_m^{-1}\Fm$, where $\NN^{-1}\!\!\approx \Am^{-1}$, and
the applications of $\NN_m^{-1}$ can be carried out cheaply.
We employ {\em two-level 
additive} and/or {\em hybrid Schwarz} techniques \cite[Chapter~2]{ToselliWidlund-DD05},
whose application in the context of discontinuous Galerkin discretization 
is very straightforward and efficient.

\subsection{Domain partition}

\begin{figure}
  \begin{center}
    \includegraphics[width=0.5\textwidth]{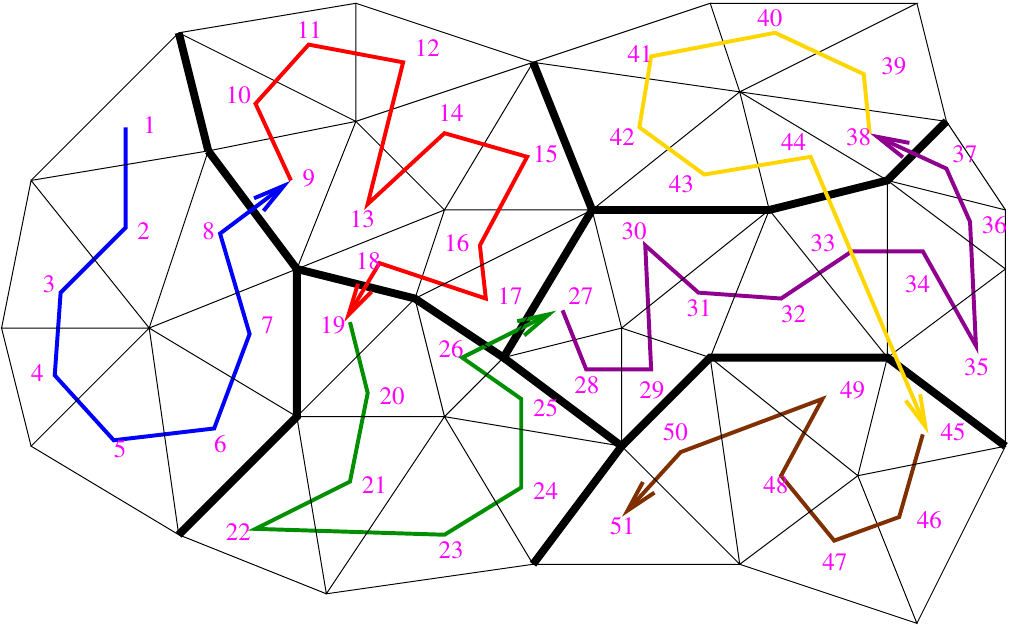}
  \end{center}
  \caption{Example of the numbering of mesh elements, bold lines correspond to the
  boundaries among subdomains $\Omi$, $i=1,\dots,\Mm$.}
  \label{fig:meshDD}
\end{figure}

Let $\Omi$ be open subdomains forming
a non-overlapping domain decomposition of the computational
domain $\Om$ at time level $i=1,\dots,m$ such that $\Omi\cap\Omj=\emptyset$ for $i\not = j$,
and $\oOm = \cup_{i=1,\dots, \Mm}\oOmi$.
The subdomains $\Omi$ are defined as a union of some (typically adjacent) elements $K\in\Thm$.
An approach to define $\Omi$, $i=1,\dots,\Mm$ is to employ graph partitioning software, e.g.,
METIS \cite{metis}, the corresponding sub-meshes are denoted as $\Thmi,\ i=1,\dots, \Mm$.
We assume that the basis functions from  $\Bm$ (cf.~\eqref{basis1}) are numbered so that 
first we number the test functions with support in $\Om_{m}^{1}$,
then the functions with support in $\Om_{m}^{2}$, etc.
The numbering is illustrated in Figure~\ref{fig:meshDD}.

\begin{figure}
  \includegraphics[width=0.330\textwidth]{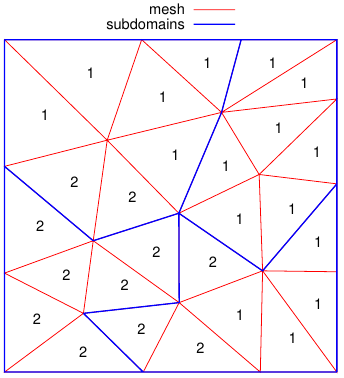}
  \includegraphics[width=0.330\textwidth]{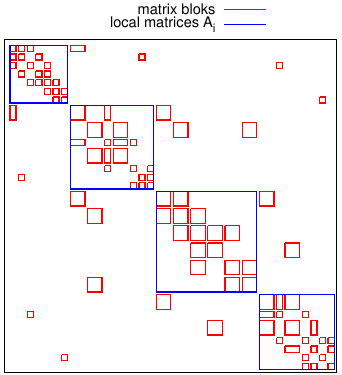}
  \includegraphics[width=0.330\textwidth]{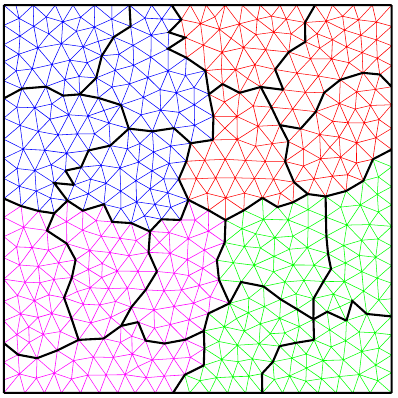}
  \caption{Example of a triangulation $\Thm$
    and the block structures of the corresponding matrix $\Am$.
    Left:  triangles with polynomial approximation degrees and the domain partition into four
    subdomains (blue lines). Center: element blocks of $\Am$ (small red boxes) and
    the subdomain-blocks  $\Ami$, $i=1,\dots,4$ (large blue boxes).
    Right: domain decomposition having four subdomains (colored) and
    the coarse grid $\THm$ with 16 coarse elements (thick lines).}
  \label{fig:blocks}
\end{figure}

In virtue of this numbering, the matrix $\Am$ from \eqref{newton1} has a natural
subdomain-block structure. 
Figure~\ref{fig:blocks} (left and center) shows an example of a triangulation
and the subdomain-block structure of the corresponding matrix $\Am$.
The left figure shows a mesh
with polynomial approximation degrees and the domain partition into four
subdomains $\Omi$, $i=1,\dots,4$. The central figure shows the block structure of
matrix $\Am$ (small red boxes).
At each block row, there is a square diagonal element-block and up to 3 off-diagonal element-blocks
corresponding to a pair of neighboring triangles.
These blocks are typically full.

Moreover, matrix $\Am$ consists of subdomain-blocks 
$\Amij$, $i,j=1,\dots,\Mm$, which correspond to terms (cf.~\eqref{newton1})
\begin{align}
  \label{Aij}
  \Amij:\quad
  \AhmL(\cdot, \bpsi_{{\bbeta}},\bpsi_{{\balpha}}),\qquad K_{\balpha}\subset\Omi,\quad  K_{\bbeta}\subset \Omj,
\end{align}
where $K_{\balpha}$ and $K_{\bbeta}$ are elements defining multi-indices
${\balpha}$ and ${\bbeta}$, respectively,
cf. Section~\ref{sec:dgm1}.
Obviously, $\Amij$ are vanishing if $\Omi$ and $\Omj$ 
do not have a common boundary
including the case when  $\Omi$ and $\Omj$ touch only in one point.
Additionally, the diagonal subdomain-blocks $\Amii$
correspond to elements from the same subdomain $\Omi$,
$i=1,\dots,\Mm$, cf. Figure~\ref{fig:blocks} (center). To shorten the notation,
we set $\Ami:=\Amii$, $i=1,\dots,\Mm$.

The space $\bSm$ (cf.~\eqref{bSm}) is naturally decomposed into local spaces on $\Omi\times I_m$,
\begin{align}
  \label{bSmi}
  \bSmi & := \left\{ \bpsi:\Omi\times I_m \to \R^n;\ \ 
  \bpsi = \bvas|_{\Omi\times(0,T)},\ \bvas\in \bSm \right\}, \qquad \ i=1,\dots,\Mm,\ \ m=0,\dots,r.
\end{align}
We denote $\Nhmi=\dim\bSmi$.
Let $\mRi\in \R^{\Nhm\times\Nhmi}$,  $i=1,\dots,\Mm$, be the matrix representation of
the restriction operator from $\bSm$ to $\bSmi$ by the identity on $\Omi$ and
by zero on $\Omj$, $j\not=i$.
Denoting $\mRiT$ the transposed operator of $\mRi$, the diagonal subdomain-blocks fulfil
\begin{align}
  \label{Ami}
  \Ami = \mRi \Am\mRiT,\qquad i=1,\dots,\Mm.
\end{align}

\subsection{Coarse subspace}

To accelerate the  transmission of information among subdomains,
a two-level method is employed.
Therefore, for each $m=1,\dots,r$, we consider a coarse mesh $\THm=\{\cK\}$, which 
consists of polygonal elements $\cK$.
We define $\cK\in\THm$ by dividing every subdomain $\Omi$ to several ``sub-subdomains'', each
of them forming one element $\cK$, see Figure~\ref{fig:blocks} (right).
Then any $\cK\in\THm$ belongs to only one subdomain $\Omi$
including the case $\cK=\oOm_i$ for some $i=1,\dots,\Mm$.
In this paper, we restrict to the case when each subdomain $\Omi$, $i=1,\dots,\Mm$, is split into
a constant number of $\sm$ coarse elements,
$\sm\in\mN$. Hence, the coarse mesh has $\#\THm=\sm\,\Mm$ elements for all
$m=1,\dots,r$.

Let $\qK= \min_{K\subset \cK} \pK$, $\ddcK:=(\qK+1)(\qK+2)/2$ for  $\cK\in\THm$,
we define the coarse finite element spaces
\begin{align}
  \label{SHp}
  \bSHm & := \left\{ \bPsi:\Om\times I_m \to \R^n;\ \ 
  \bPsi(x,t)|_{\cK\times I_m} = \sum\nolimits_{j=0}^q t^j \bVas(x),\
  \bVas\in [P_{\qK}(\cK)]^n,\  \cK\in\THm\right\},
  \quad \ m=0,\dots,r,
\end{align}
cf.~\eqref{bSm}.
Similarly as in \eqref{basis1}, we define the basis
$\BHm$ of spaces $\bSHm$ by replacing the
symbol $K$ by $\cK$ and $h$ by $H$, i.e.,
\begin{align}
  \label{basisH}
  \BHm&:=\left\{\bPsi_{\balpha},\ \bPsi_{\balpha}:= \bPsi_{\cK,k,i}^{j},\  {\balpha}=\{\cK,k,i,j\}\ \
  \forall {\balpha} \in \BBHm\right\},
\end{align}
where $\BBHm$ is the set of all multi-indices
${\balpha}=\{\cK,k,i,j\}:\ j=0,\dots,q,\  i=1,\dots,\ddcK,\ k=1,\dots,n,\ \cK\in\THm$.

Obviously,  $\bSHm\subset \bSm$ for all $m=1,\dots, r$.
Therefore, any basis function from $\BHm$, can be expressed as a linear combination
of basis function from $\Bm$, cf.~\eqref{basis1}. The coefficients of these linear combinations
defines the {\em restriction matrix} $\mRN\in\R^{\NHm\times\Nhm}$ where
$\NHm$ is the dimension of $\bSHm$ given analogously as \eqref{Nhm}. Let $\mRNT$ be its transpose
representing the prolongation operator. Then we set the coarse matrix
\begin{align}
  \label{A0}
  \AmN := \mRN \Am \mRNT \in \R^{\NHm\times\NHm},\qquad  m=1,\dots,r.
\end{align}

\begin{remark}
  \label{rem:A0}
  In virtue of \eqref{newton1} and the definition of $\mRN$,
  the matrix $\AmN$ from \eqref{A0} fulfils
  \begin{align}
    \label{newtonA0}
    \AmN(\bW) = \left\{ (\AmN(\bW))_{{\balpha},{\bbeta}} \right\}_{{\balpha},{\bbeta}\in\BBHm},
    \quad
    (\AmN(\bW))_{{\balpha},{\bbeta}} = \AhmL(\bbw, \bpsi_{{\bbeta}},\bpsi_{{\balpha}}),\quad {\balpha},{\bbeta}\in\BBHm.
  \end{align}
  Hence, $\AmN$ can be evaluated directly from \eqref{Ah} inserting the basis function from
  $\BHm$ to the second and third arguments of $\AhmL$.
  Although both \eqref{A0} and \eqref{newtonA0} are equivalent,
  the approach \eqref{newtonA0} is less robust in the finite precision arithmetic
  when solving the Navier-Stokes equations.
  On the other hand,
  we have not observed this effect for scalar problems ($n=1$ in \eqref{eq:NS}).
\end{remark}

\subsection{Additive and hybrid preconditioners}
\label{sec:precond}

We introduce the additive and hybrid two-level Schwarz preconditioners
which are suitable for parallelization. We present only the final formulas, for
details we refer, e.g., to\cite[Chapter~2]{ToselliWidlund-DD05}.
The {\em additive preconditioned operator}
is additive in the local problems
as well as the global coarse one, and it reads
\begin{align}
  \label{ASM2}
  \mPadd  :=  \sum\nolimits_{i=0}^{\Mm}\mPPi,\qquad
  \mbox{where} \quad \mPPi:=  \mRiT (\Ami)^{-1}\mRi \Am,\quad i=0,\dots,\Mm.
\end{align}
Since $\mPadd  = \mNadd^{-1} \Am$, the corresponding preconditioner has the form
\begin{align}
  \label{ASM1}
  \mNadd^{-1} =  \sum\nolimits_{i=0}^{\Mm}  \mRiT (\Ami)^{-1}\mRi.
\end{align}

%

On the other hand, the hybrid preconditioner is additive in the local components but
multiplicative with respect to the levels.
The {\em hybrid preconditioned operator} reads
\begin{align}
  \label{hyb2}
  \mPhy &:= \mI - \big(\mI - \mPPN \big)\big(\mI - \sum\nolimits_{i=1}^N \mPPi \big),
\end{align}
where $\mPPi$ are given in \eqref{ASM2}. Since $\mPhy :=  \mNhy^{-1} \Am$,
the hybrid preconditioner can be written in the form
\begin{align}
  \label{hyb1}
  \mNhy^{-1}  
  & =  \sum\nolimits_{i=1}^N \mRiT (\Ami)^{-1} \mRi
  + \mRNT (\AmN)^{-1}\mRN \Big(\mI - \Am \sum\nolimits_{i=1}^N \mRiT (\Ami)^{-1} \mRi\Big).
\end{align}
We note that the hybrid preconditioned operator \eqref{hyb2} is not symmetric. Therefore,
a modification for symmetric problems is required. However, the systems treated in our paper are
not symmetric so the form \eqref{hyb2} is applicable.

\begin{algorithm}[t]
  \caption{Application of the preconditioner $\mNadd^{-1}$ from \eqref{ASM1}:
    $\bku \leftarrow \mNadd^{-1}\bkx$}
  \label{alg:ASM}
  \begin{spacing}{1.15}
    \begin{algorithmic}[1]
      \STATE \textbf{input} matrices $\Am$, $\Ami$, $\mRi$, $i=0,\dots,N$, vector $\bkx$
      \STATE set $\bkx_i := \mRi \bkx$ and solve $\Ami \bky_i = \bkx_i$ for $i=0,\dots,N$
      \STATE \textbf{output} vector  $\bku := \sum_{i=0}^N \mRi^T \bky_i$
    \end{algorithmic}
  \end{spacing}
\end{algorithm}

\begin{algorithm}[t]
  \caption{Application of the preconditioner $\mNhy^{-1}$ from \eqref{hyb1}:
    $\bku \leftarrow \mNhy^{-1}\bkx$}
  \label{alg:ASM2}
  \begin{spacing}{1.15}
    \begin{algorithmic}[1]
      \STATE \textbf{input} matrices $\Am$, $\Ami$, $\mRi$, $i=0,\dots,N$, vector $\bkx$
      \STATE set $\bkx_i := \mRi \bkx$ and  solve $\Ami \bky_i = \bkx_i$ for $i=1,\dots,N$
      \STATE $\bky := \sum_{i=1}^N \mRi^T \bky_i$
      \STATE $\bkz := \bkx -\Am \bky $ \label{asm2_m1}
      \STATE $\bkx_0 := \mRN \bkz$ and solve $\AmN \bky_0 = \bkx_0$ \label{asm2_c1}
      \STATE \textbf{output} vector $\bku := \mRN^T \bky_0 + \bky$
    \end{algorithmic}
  \end{spacing}
\end{algorithm}

Algorithms~\ref{alg:ASM} and \ref{alg:ASM2}
describe
the applications of preconditioners $\mNadd^{-1}$ and $\mNhy^{-1}$
introduced in
\eqref{ASM1} 
and \eqref{hyb1}. 
While Algorithm~\ref{alg:ASM} allows  to solve the local fine problems
(with $\Ami$, $i=1,\dots,N$) together  with the global coarse problem (with $\AmN$)
in parallel,
Algorithm~\ref{alg:ASM2} requires solving the local fine problems first
and then solving the global coarse problem. Additionally, one multiplication by
$\Am$ has to be performed.
The multiplication by $\Am$ is typically cheaper than
the solution of the local or coarse problems. Moreover, if the global coarse problem is smaller
than a local one, then the application of the preconditioner $\mNhy^{-1}$
exhibits only  a small increase of the computational time in comparison to
the preconditioner $\mNadd^{-1}$.
Additionally, the hybrid preconditioner is typically
more efficient, and a smaller number of GMRES iterations
has to be carried out to achieve the given tolerance.


The application of the two-level preconditioners in Algorithms~\ref{alg:ASM} and \ref{alg:ASM2}
exhibits the solution of the local (fine) systems and the global (coarse) one,
\begin{align}
  \label{flops0}
  \Ami \bky_i = \bkx_i,\quad i=1,\dots,\Mm,\qquad \mbox{and} \qquad
  \AmN \bky_0 = \bkx_0, \qquad \mbox{respectively}.
\end{align}
We solve these systems directly by the MUMPS library \cite{MUMPS,MUMPS1,MUMPS2},
which consists of two steps:
\begin{enumerate}[label=(S{\alph*})]
\item the {\em factorization} of the system, which is carried out
  after the evaluation of $\Am$, \label{P1}
\item the {\em substitution}
  of the solution which is carried out in each GMRES iteration.
  \label{P2}
\end{enumerate}
We note that the term ``substitution'' means the evaluation of the solution using
the factorization of the matrix and it is called ``assembling'' in MUMPS. However, in the
finite element community, the assembling usually means the creation of the matrix itself,
so we use the term substitution.

\subsection{Problem formulation}
\label{sec:form}
The numerical method described in Section~\ref{sec:DGM} requires the solution of a sequence of
linear algebraic systems 
\begin{align}
  \label{addm1}
  \Am(\WmlM) \dl = \Fm(\WmlM),\qquad \ell=1,2\dots, \quad m=1,\dots,r,
\end{align}
where $\ell$ denotes the Newton iterates and $m$ the time level, cf.~\eqref{newton2b}.
We solve \eqref{addm1} by GMRES with two-level preconditioners described
in Section~\ref{sec:precond}. Finally, the GMRES solver is stopped when the norm of the
preconditioned residual decreases by a fixed factor $\CL\ll1$.
Hence, the algebraic error is controlled by \eqref{CA} only.

These preconditioners have to solve
$\Mm$ local problems defined on subdomains $\Omi$, $i=1,\dots,\Mm$, and the global problem
on a coarse mesh having $\#\THm = \sm\,\Mm$ elements, where $\sm\ge 1$ is
a chosen integer.
We recall that  we consider the case when each $\Omi$, $i=1,\dots,\Mm$,
is divided into $\sm$ coarse elements.
The fundamental question studied in this paper is the following.

\begin{problem}
  \label{prob:main}
  Let $\bSm$, $m=1,\dots,r$ be, the sequence of finite element spaces \eqref{bSm}
  generated by {\hpAMA}
  (cf. Section~\ref{sec:AMA}) used for the solution of \eqref{eq:NS}
  by STDGM (cf. Section~\ref{sec:DGM}).
  How to choose $\Mm$ and $\sm$, $m=1,\dots,r$, such that the solution of the arising sequence
  of linear systems \eqref{addm1} by GMRES with 
  two-level preconditioners requires the shortest possible
  wall-clock time in a parallel computation?
\end{problem}

Obviously, larger $\Mm$ admits to use more computer cores,
but the GMRES solver requires more iterations. On the other hand,
fines coarse mesh $\THm$ increases the speed of convergence,
but the corresponding coarse problem is more expensive to solve.
Moreover,
the real wall-clock time depends on many factors.
Apart of the numerical scheme itself, also on
its implementation, the speed of communication among the cores of the computer, their type, etc.
Therefore, we introduce a simplified computational cost model which 
takes into account the number of floating point operations, the speed of computation
(performance), and
the wall-clock time of the communication among the cores of the computer.

\section{Computational costs}
\label{sec:costs}

In this section, we introduce the simplified model which will measure
the computational costs of the preconditioners which exhibits the most time-consuming
part of the whole process. Based on this model, we propose an adaptive algorithm
giving an approximation of Problem~\ref{prob:main}.

\subsection{Quantities defining the computational costs}
To define the computational cost model, we introduce the following assumptions.
\begin{enumerate}[label=(A{\arabic*})]
\item  \label{ass1}
  The applications of the two-level preconditioner,
  i.e., the solution of systems \eqref{flops0},
  is the dominant part of the computational process. The other parts, such as
  the evaluation of entries of matrix $\Am$ in~\eqref{newton1}  and vector $\Fm$
  in~\eqref{newton0},
  multiplication of a vector by matrix $\Am$ in GMRES, the computational overheads, etc.,
  can be neglected.
\item Enough computer cores are available such that each local system in \eqref{flops0}
  can be solved in parallel, and also each $\Ami$, $i=0,\dots,\Mm$
  is stored only at a separate core. Therefore, the factorization of the local systems
  and the global coarse one can be carried out simultaneously in parallel.
  \label{ass2}
\item All vectors appearing in the computations are stored in copies at each used computer core.
  Therefore, in virtue of the previous item, the communication among the cores is represented just
  by the output of~\eqref{flops0} sent from each core and received by every core.
  \label{ass3}
\item Only the following communications among the cores are considered:
  \begin{enumerate}
  \item The solution of the local systems \eqref{flops0}
    for $i=1,\dots,\Mm$ gives the local vectors $\bky_i$ which have to be joined together
    and spread among all cores. We denote this communication $\gathr$ 
    and use the subroutine  {\it MPI\_Allgatherv} in MPI \cite{MPI50}.
  \item The solution of the coarse system \eqref{flops0}
    for $i=0$ gives the local vector $\bky_0$ which has to be spread among all cores.
    We denote this communication $\bcast$ 
    and use the subroutine  {\it MPI\_Bcast} in MPI \cite{MPI50}.
  \end{enumerate}
  \label{ass4}
\end{enumerate}
We are aware that the previous assumptions exhibit a strong simplification, yet
they allow us to define a relatively simple computational cost model which
can provide useful information.
The previous assumptions can be relaxed by incorporating the corresponding
items in a more complex model.

\begin{remark}
It is possible to consider an alternative implementation, in which the coarse problem
is stored in copies at each core. Then  $\bcast$ can be avoided but the  coarse problem
can not be factorized simultaneously with the local ones. Consequently,
the  subsequent computational model would have to be modified. 
\end{remark}

The goal of the computational cost model is not only to measure the computational costs
but also to predict them for a priori unknown domain decomposition.
Therefore, we consider 
\begin{enumerate}[label=(T{\arabic*})]
\item the number of {\em floating point operations} (\flops) per each
  task of the algorithm  and each core, \label{term1}
\item the {\em speed of the computations} (= performance) per each task of the algorithm
  and each core,
  \label{term2}
\item the wall-clock time of {\em communication operations} among the cores.
  \label{term3}
\end{enumerate}
We note that the number of {\flops} necessary to solve \eqref{flops0} are provided
by the MUMPS library, cf. Section~\ref{sec:flops}.

\subsection{Simplified model of the algorithm on parallel computer}

\begin{figure} [t]
  \begin{center}
    \includegraphics[width=0.85\textwidth]{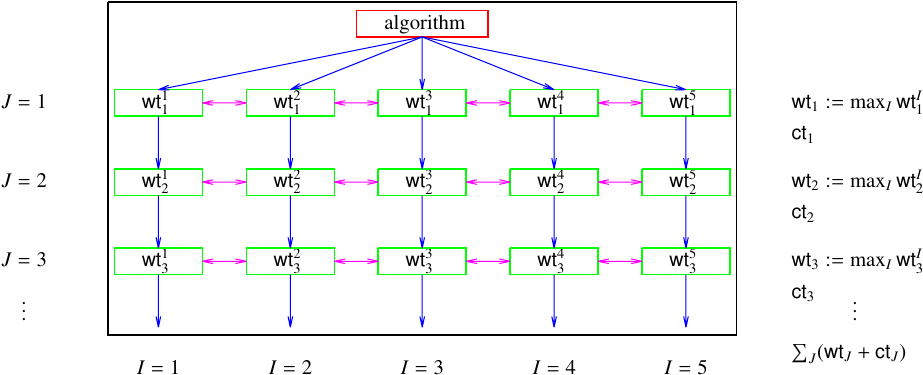}
\end{center}
  \caption{Computational cost  model using 5 cores  $I=1,\dots,5$,
    parallel computations (green boxes)
    followed by communication among cores (magenta arrows), and subsequent steps of the
    algorithm (blue arrows).}
  \label{fig:computer}
\end{figure}

For the setting of the computational costs, we use an algorithm implementation sketched in
Figure~\ref{fig:computer}. The algorithm is split into several task levels $J=1,2,\dots$,
and each task level is split into $M$ tasks, one task for one core $I$, $I=1,\dots,M$.
The performance of each task on each core is independent of the performance of
other tasks and cores. We denote by
$\task_J^I$ the task on task level $J=1,2,\dots$ and the core $I=1,\dots,M$.
When all tasks $\task_J^I$, $I=1,\dots,M$ from one task level $J$ are finished, the cores
are synchronized, 
the particular outputs are transmitted among the cores, and we proceed to the next task level $J+1$.
Let $\Tfl_J^I$ denote 
the wall-clock time 
to carry out task $\task_J^I$,
$I=1,\dots,M$, $J=1,2,\dots$,  then we define the wall-clock time 
for task level $J$ by
\begin{align}
  \label{flopsB1}
  \Tfl_J := \max_{I=1,\dots,M} \Tfl_J^I,\qquad J=1,2,\dots.
\end{align}
By $\Cfl_J$ we denote the wall-clock time necessary to transmit data from level $J$ to level
$J+1$. Consequently the total 
wall-clock time
of the algorithm
is given by
\begin{align}
  \label{flopsB2}
  \Tfl := \sum\nolimits_{J=1,2,\dots} (\Tfl_J+\Cfl_J).
\end{align}
Note that that the measure $\Tfl_J$ can be seen as the wall-clock time  needed by one core with
the largest amount of work (longest time) at level $J$. 
In the following, we study the particular task levels of the algorithm in terms
of \ref{term1}--\ref{term3} under Assumptions~\ref{ass1}--\ref{ass4}.
Therefore, we consider only the factorization and substitution of the local and
coarse global systems \eqref{flops0}. 

\subsection{Computational cost  model}
\label{sec:model0}
First, we introduce the following notation which relates to the discretization and
the computational process introduced in Sections~\ref{sec:DGM} and \ref{sec:DDM}.
Let $m=1,\dots,r$ be arbitrary but fixed index of the time layer.
\begin{itemize}
  \itemsep0pt 
\item $\#\Thm$ -- number of elements of mesh $\Thm$,
\item $\Nhm$ ($=\dim\bSm$) --  the size  of matrix $\Am$,
\item $\Mm$ -- the number of subdomains of domain decomposition,
\item $\sm$ -- the number of coarse elements in each subdomain, the number of coarse elements
  is $\#\THm = \sm\, \Mm$,
\item $\iterN$ -- the number of Newton (nonlinear) iterations in time step $m$,
\item $\iterL$ -- the number of GMRES (linear) iterations over all
  nonlinear iterations $\iterN$ in time step $m$,
\item $\callC$ -- the number of evaluations of matrix $\Am$ from \eqref{newton1}
  in time step $m$. 
\end{itemize}
We note that $\callC \leq \iterN$ and the equality is valid if the matrix $\Am$ is updated
at each Newton iteration. However, in virtue of Remark~\ref{rem:refresh}, the value
$\callC$ is smaller in practice.
In the following paragraphs,
we discuss the particular part of the computation in terms \ref{term1}--\ref{term3}
in one time step $m$.

\subsubsection{Floating point operations}
\label{sec:flops}

The solution of systems
\eqref{flops0} by MUMPS  consists of the factorization of the matrices (step \ref{P1})
and the substitution of the solution (step \ref{P2}).
While step \ref{P1} is performed only after the refreshing of $\Am$ (i.e., $\callC$ times),
the substitution is carried out at each GMRES iterations.

Let $\flfaci$ and $\flassi$ denote the number of floating point operations
for factorization of $\Ami$ and substitution of the corresponding solutions,
respectively, $i=0,\dots,\Mm$.
The factorization of $\Ami$ can be carried out independently for
each $i=0,\dots,\Mm$, so in virtue of \eqref{flopsB1},
the maximum of the number of floating point operations per core is 
\begin{align}
  \label{flopsF1}
 \FFfac := \max\nolimits_{i=0,\dots,\Mm} \flfaci.
\end{align}
Note that index 0 is included and stands for the coarse problem.
On the other hand, the substitution of the solutions of the local systems and the global coarse one
\eqref{flops0}
is executed in parallel only for the additive preconditioner $\mNadd^{-1}$ \eqref{ASM1}.
For the hybrid preconditioner  $\mNhy^{-1}$ \eqref{hyb1}, we solve the local systems in parallel
and then the global system sequentially. Hence, the corresponding
maximum of the number of floating point operations per core for substitution is
\begin{align}
  \label{flopsF2}
  \FFass :=
  \begin{cases}
    \max_{i=0,\dots,\Mm} \flassi & \mbox{ for } \mNadd^{-1},\\
    \max_{i=1,\dots,\Mm} \flassi  + \flassN & \mbox{ for } \mNhy^{-1},\\
  \end{cases}
\end{align}
in each GMRES iteration. Therefore, in virtue of \eqref{flopsB1}--\eqref{flopsB2},
the total number of $\fl$ per one node
is given by
\begin{align}
  \label{flopsF3}
  \flm = \callC\, \FFfac + \iterL \, \FFass,\qquad m=1,\dots,r.
\end{align}
We note that the substitution at each time step is carried out
$\iterL +  \iterN$ times since the preconditioner has
to be applied also to the right-hand side of \eqref{addm1} in each Newton iterate.
However, the value $\iterL$ is typically several times larger than  $\iterN$, hence the
latter can be omitted in \eqref{flopsF3}.

\begin{figure} [t]
  \begin{center}
    \includegraphics[width=0.49\textwidth]{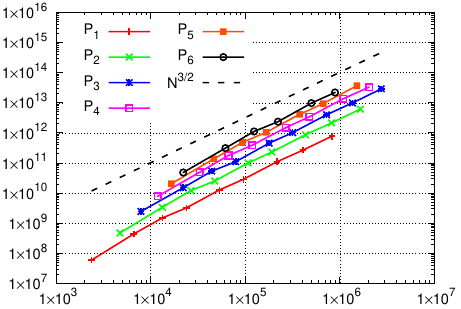}   
    \includegraphics[width=0.49\textwidth]{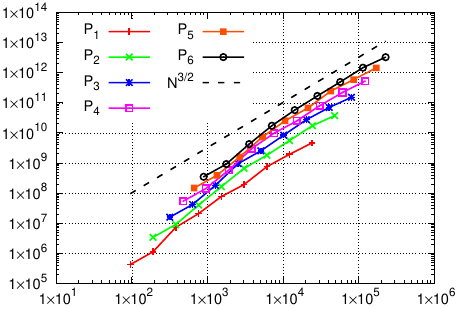}
\end{center}
  \caption{Number of the floating point operations for the matrix factorization:
    dependence of $\flfaci$ on the  system size $\Nhmi$ 
    corresponding to triangular $\Thmi$ (left) and polygonal $\THm$(right) grids.}
  \label{fig:flops0}
\end{figure}
\begin{figure} [t]
  \begin{center}
    \includegraphics[width=0.48\textwidth]{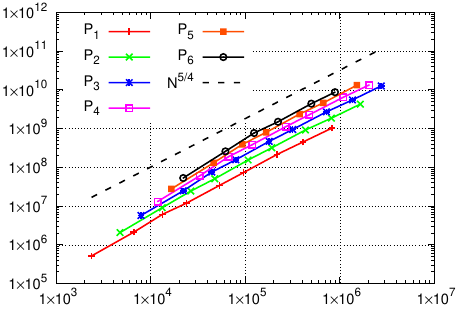}
    \hspace{0.03\textwidth}
    \includegraphics[width=0.48\textwidth]{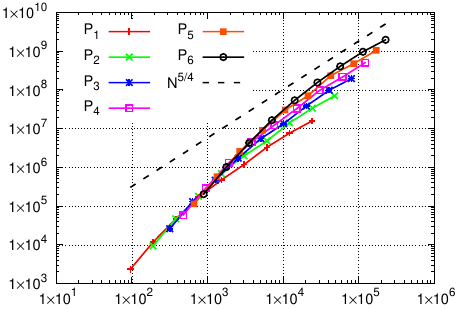}
\end{center}
  \caption{Number of the floating point operations for the solution substitution:
    dependence of $\flassi$  on the system size $\Nhmi$ corresponding to
    triangular $\Thmi$ (left) and $\THm$ polygonal (right) grids.}
  \label{fig:flops2}
\end{figure}

To establish the computational cost model, we are interested in
the dependence of $\flfaci$ and $\flassi$ on the size of the
systems $\Nhmi$, $i=0,\dots,\Mm$. These values depend on the sparsity and shape
of the matrices $\Ami$ and therefore on meshes $\Thm$, $\THm$, polynomial degrees, etc.
We carried out a series of numerical experiments when the Navier-Stokes equations \eqref{eq:NS}
were discretized by STDG method~\eqref{STDGM} on several quasi-uniform meshes $\Thm$ and
coarse meshes $\THm$
with fixed temporal degree $q=1$ and uniform spatial degrees $\pK=1,\dots,6$, $K\in\Thm$.
The resulting algebraic systems $\Am$ 
and $\AmN$ were factorized and substituted by MUMPS
(cf.~\ref{P1}--\ref{P2}) which provides the corresponding number of $\fl$
(variables {\it mumps\_par\%RINFOG(1:2)} in \cite[User's guide]{MUMPS}).
We consider the coarse meshes $\THm$ separately since they are polygonal and
their sparsity  differs from the sparsity (fine or local) meshes of triangular elements.

Numerical results are given in Figures~\ref{fig:flops0} and \ref{fig:flops2},
they  show  the dependence of the values $\flfaci$ and $\flassi$  on the system size,
respectively,
in logarithmic scale. We easily deduce an asymptotic exponential dependence
\begin{align}
  \label{addm5}
  \flfaci(N) \approx C \, N^\mu,\qquad \flassi(N) \approx c\, N^\nu, \qquad
  i=1,\dots,\Mm,
\end{align}
where $c$, $C$, $\mu$ and $\nu$ are positive parameters. More detailed inspection
indicates the values $\mu\approx 3/2$ and $\nu\approx 5/4$ for all polynomial degrees $p$
whereas the constants $c$ and $C$ are $p$-dependent. It is caused by the increasing sparsity
of matrices for increasing $p$.
The empirical formulas \eqref{addm5} are used later for the prediction
of the number of floating point operations.


\subsubsection{Speed of the computation}
\label{sec:speed}

\begin{figure} [t]
  \begin{center}
    \includegraphics[width=0.49\textwidth]{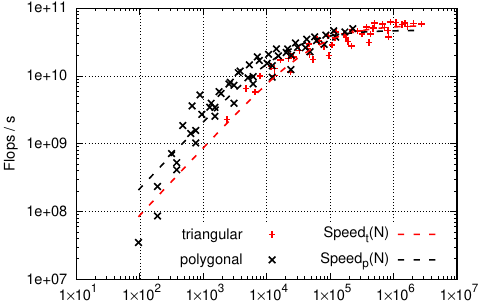}   
    \includegraphics[width=0.49\textwidth]{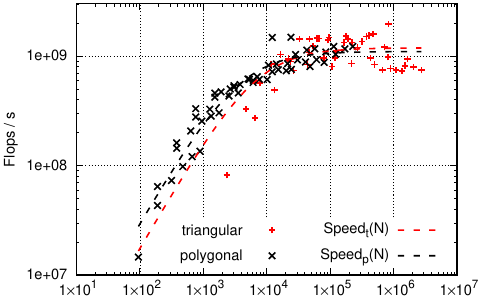}
\end{center}
  \caption{Performance of the matrix factorization and substitution:
    dependence of $\flfaci\,/\, s$ (left) and $\flassi\, /\, s$ (right)
    on the  system size $\Nhmi$ for triangular $\Thmi$ (red +) and polygonal $\THm$
    (black $\times$). Dashed lines are fitted functions $\speed(N)$ in the form of \eqref{addm6} for
    both types of meshes and factorization/substitution.}
  \label{fig:speed0}
\end{figure}

Although the number of $\fl$ discussed in Section~\ref{sec:flops} gives reasonable information
about the computational costs, it does not fully reflect the wall-clock time which is
the desired interest in practical computations.
The wall-clock time does not scale linearly with the number of $\fl$
typically due to limits to speed of memory access.
Therefore, we measure the wall-clock time  for the factorization and substitution
of the systems in all numerical experiments from Section~\ref{sec:flops}.
We use the function {\it MPI\_Wtime} in the MPI \cite{MPI50} which returns
the wall-clock time in seconds ($s$). We denote the wall-clock times for the factorization and
substitution of~\eqref{flops0}  by $\Wtfaci$ and $\Wtassi$, $i=0,\dots,\Mm$, respectively.
Then the speed of computation of factorization and substitution reads
\begin{align}
  \label{addmS1}
  \spfaci := \frac{\flfaci}{\Wtfaci}  \qquad \mbox{and} \qquad
  \spassi := \frac{\flassi}{\Wtassi}, \qquad i=0,\dots,\Mm,
\end{align}
respectively. 
Figure~\ref{fig:speed0} shows the dependence of $\spfaci$ and $\spassi$
on the size of the system $\Nhmi$ in units $\fl/s$.
The numerical experiments were carried out using
cluster of nodes with 4x Intel Xeon Gold 6240 CPU \@ 2.60GHz, 512GB RAM,
interconnected by InfiniBand 100 Gb.

Figure~\ref{fig:speed0} shows that the performance is increasing for increasing size of the
system and it is approaching a limit of a maximal performance.
Therefore, we estimate a simple empirical dependence
\begin{align}
  \label{addm6}
  \speed(N) \approx  \frac{a\, N}{ b + N},\qquad a,b > 0,
\end{align}
where $a$ is the maximal speed for $N\to\infty$. The parameters $a$ and $b$ can be found by
fitting the function in the form \eqref{addm6} with the measured data.
We have found them 
iteratively by a nonlinear least square technique, and
the corresponding values of $a$ and $b$ from \eqref{addm6} are given in Table~\ref{tab:fitting}.
The fitted functions are shown also in Figure~\ref{fig:speed0} by dashed lines.

\begin{table}
  \caption{Performance of the matrix factorization and substitution, the values of the
    fitted parameters $a$ and $b$ from \eqref{addm6}.}
   \label{tab:fitting}
  \begin{center}
    \setlength{\tabcolsep}{12pt}
    \input{Figs_flops/speed.tex}
  \end{center}
\end{table}

The measuring of the wall-clock time is influenced by several perturbations
caused by the system noise, operating system scheduling, etc.
These tiny interruptions lead to slight variations in the wall-clock time, up to several percents.
Moreover, we observed that the fitting of $\speed$ according to formula \eqref{addm6} by the
nonlinear least-squares technique is sensitive to data perturbations and initial approximation.
Therefore, in practical applications, we replace~\eqref{addm6} by an affine
approximation, particularly
\begin{align}
  \label{addm7}
  \speed(N) \approx  \min(\Smax, {a\, N} +  b),
\end{align}
where $a$ and $b$ are the parameters given by a linear least-squares method and
$\Smax$ is the maximal speed. In practical computation, we set $\Smax$
according to formula \eqref{fit2} below.
Formula \eqref{addm7} mimics dependence \eqref{addm6} locally, and in our experience,
it is less sensitive
to data perturbations. 

Finally, we recall that from \eqref{addmS1}, we can estimate time as
\begin{align}
  \label{addm8a}
   \Wtfaci = \frac{\flfaci}{\spfaci},\qquad
   \Wtassi = \frac{\flassi}{\spassi}, \qquad i=0,\dots,\Mm,
\end{align}
and therefore, using \eqref{flopsF3}
the wall-clock time for the application of the preconditioners for all GMRES iterations
within time step $m$ is estimated as
\begin{align}
  \label{addm8}
  \Wtimem &=
    \callC\, \max_{i=0,\dots,\Mm} \Wtfaci  + \iterL \, \max_{i=0,\dots,\Mm} \Wtassi
    \hspace{23mm} \mbox{ for } \mNadd^{-1},\\[2mm]
    \Wtimem &=
    \callC\, \max_{i=0,\dots,\Mm} \Wtfaci  + \iterL \Big( \max_{i=1,\dots,\Mm} \Wtassi +  \WtassN\Big)
    \qquad  \mbox{ for } \mNhy^{-1}.
    \notag
\end{align}

\subsubsection{Communication operations} 
\label{sec:comm}
In virtue of Assumption~\ref{ass4}, we consider only the communications
$\gathr$ and $\bcast$ which are carried out always after the solution of systems \eqref{flops0},
hence once in each GMRES iteration. These communication costs
have two fundamental parts: the {\em latency}, which
is the fixed cost of starting a communication, independent of
the size of transmitted data, and
the {\em bandwidth} which is the rate at which data can be
transferred once the communication has started.
According to, e.g., \cite{ChanALL_CCRE07,KangALL_PC19}, the lower bounds 
of the latency and bandwidth of $\gathr$ and $\bcast$ are given by  
\begin{align}
  \label{comm1}
  \mbox{latency} \approx\alpha\,\log_2(P), \qquad  \qquad
  \mbox{bandwidth} \approx \beta\, \frac{P-1}{P}\, N, 
\end{align}
where $P$ is the number of cores, and $N$ is the length of transmitted data.
The parameters $\alpha$ and $\beta$ are positive constants and they are hardware and
software dependent.

To illustrate the dependence of $\gathr$ and $\bcast$  on the number of processors $P$ and
the size of transmitted data $N$, we carried out several numerical experiments measuring
the wall-clock time of both communications.
We denote by  $\TgathrI$  and $\TbcastI$ the wall-clock time of $\gathr$ and $\bcast$
of one transmission given by Assumption~\ref{ass4}, respectively.
To obtain observations close to a real computation, 
we applied the numerical scheme presented
in Sections~\ref{sec:DGM} and \ref{sec:DDM} for unstructured quasi-uniform
meshes $\Th^i$, $i=1,\dots,4$
having 996, 2248, 3996, and 9074
elements. We considered 80 degrees of freedom per element, hence the size of the
corresponding global fine systems was $N=79\,680$, 179\,840, 319\,680, and 725\,920.
For each mesh, we consider $M=8$, 16, 32, and 64 subdomains, and each computation
was executed with $P=M$ cores. The number of coarse elements was chosen as $2^{i-1}\,M$
and then the size of small (global coarse) system was $N^0=80\cdot 2^{i-1}\,M$
for all meshes  $\Th^i$, $i=1,\dots,4$.
We note that the size of the local vectors (each local vector is associated to one core)
is not the same, hence we use {\it MPI\_Allgatherv} function.
The cluster setting and the measurement of the wall-clock time was the same as
in Section~\ref{sec:speed}.

\begin{figure} [t]
  \begin{center}
    \hspace{0.05\textwidth}
    $\Tbcast$(using \ {\it MPI\_Bcast})  \hspace{0.22\textwidth}
    $\Tgathr$ (using \ {\it MPI\_Allgatherv})

    \includegraphics[width=0.49\textwidth]{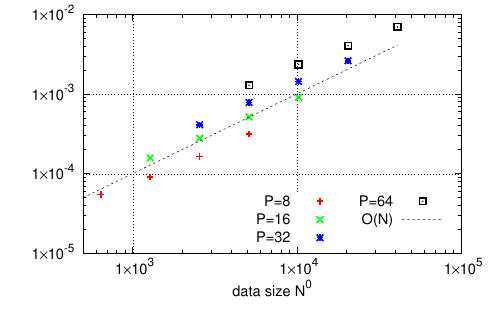}   
    \includegraphics[width=0.49\textwidth]{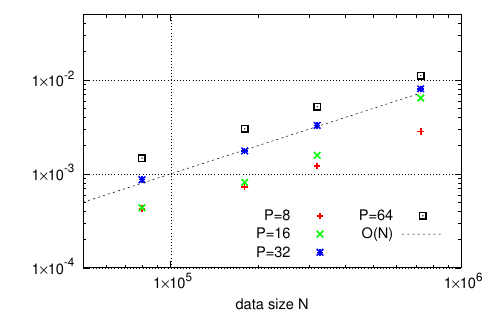}

    \vspace{2mm}
    
    \includegraphics[width=0.49\textwidth]{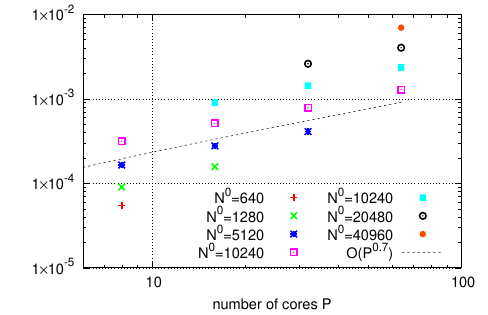}   
    \includegraphics[width=0.49\textwidth]{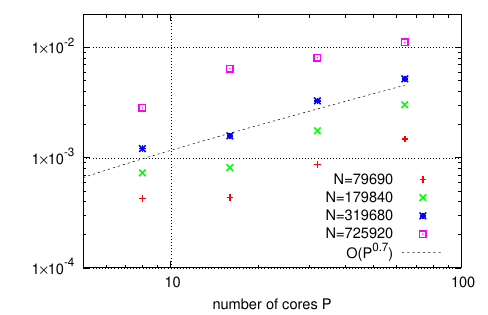}
\end{center}
  \caption{Communication wall-clock time for one
    $\bcast$ (left columns) and $\gathr$ (right columns),
    the dependence on the size of communicated vectors (top row) and the number of cores
  (bottom row).}
  \label{fig:commP}
\end{figure}

In order to avoid fluctuation in the measurement of the communication time, we carried out
at least one thousand of GMRES iterations and present averaged times.
The results are given in Figure~\ref{fig:commP}, where we show the dependence of the 
(averaged)
wall-clock time  of one  $\bcast$ and one $\gathr$ communication
on the size of transmitted vector ($N^0$ for $\bcast$, $N$ for $\gathr$)
and on the number of cores $P$.

We observe the linear dependence $O(N)$ of the wall-clock times 
for both communications $\gathr$ and $\bcast$ on the size of the transmitted vector
which is in agreement
with the bandwidth \eqref{comm1} (for large $P$ when $(P-1)/P\approx 1$).
On the other hand, the dependence of  both communications $\gathr$ and $\bcast$ on the
number of cores is clearly sub-linear but not quite logarithmic.
We found empirically the dependence $O(P^{0.7})$, which is depicted
in Figure~\ref{fig:commP} too.

\begin{figure} [t]
  \begin{center}
    \hspace{0.05\textwidth}
    $\Tbcast$(using \ {\it MPI\_Bcast})  \hspace{0.22\textwidth}
    $\Tgathr$ (using \ {\it MPI\_Allgatherv})

    \includegraphics[width=0.49\textwidth]{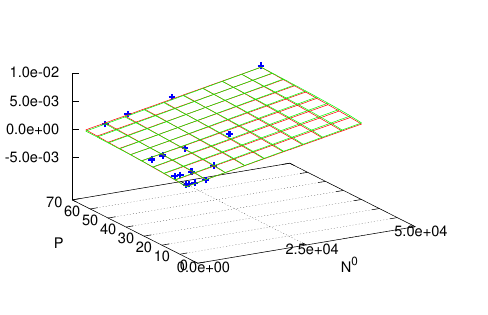}   
    \includegraphics[width=0.49\textwidth]{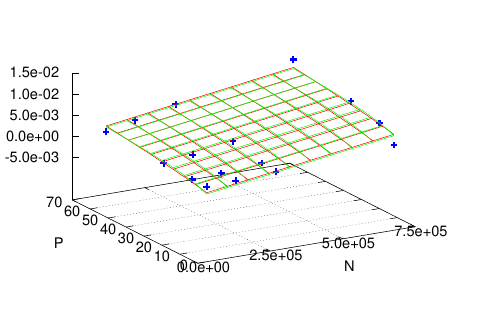}
\end{center}
  \caption{Communication wall-clock time for one
    $\bcast$ (left columns) and $\gathr$ (right columns),
    the dependence on the size of communicated vectors $N$ and the number of cores $P$
    (blue nodes), fitted function using \eqref{comm2} (green) and \eqref{comm3} (red).}
  \label{fig:comm3D}
\end{figure}

Therefore, in virtue of \eqref{comm1} and the numerical experiments,
two variants of the communication cost model are possible:
\begin{subequations}
  \label{comm0}
  \begin{align}
    \label{comm2}
    & \TgathrI(N,P) \approx \alpha_g \log_2{P} + \beta_g N + \gamma_g, \hspace{11mm}
    \TbcastI(N,P) \approx \alpha_b  \log_2{P}  + \beta_b N + \gamma_b, \\
    \label{comm3}
    & \TgathrI(N,P) \approx \alpha_g (P-1)^{0.7} + \beta_g N+ \gamma_g, \qquad
    \TbcastI(N,P)   \approx \alpha_b (P-1)^{0.7} + \beta_b N+ \gamma_b,
  \end{align}
\end{subequations}
where $\TgathrI$ and $\TbcastI$ are the wall-clock times of $\gathr$ and $\bcast$
for one call of {\it MPI\_Allgatherv} and {\it MPI\_Bcast},  respectively.
To observe the difference between \eqref{comm2} and \eqref{comm3},
we found the parameters  $\alpha_g$,  $\beta_g$,  $\gamma_g$,  $\alpha_b$, $\beta_b$, and $\gamma_b$
for the data from  Figure~\ref{fig:commP} by the least square technique. The resulting
fitted functions are shown in  Figure~\ref{fig:comm3D}, where
all data from  Figure~\ref{fig:commP} (blue nodes) are fitted by analytical functions
\eqref{comm2} and \eqref{comm3}. We observe a close approximation of the data.
In the following, we use fitting   \eqref{comm2} since it has better
theoretical background and may work better also beyond one node.

Finally, we conclude that the communications $\bcast$ and $\gathr$
are carried out once in each GMRES iteration, hence
the wall-clock time for the communication in $m$-th time step can be approximated as
\begin{align}
  \label{comm4}
  \Tgathr(N,P) = \iterL \, \TgathrI, \qquad
  \Tbcast(N,P) =  \iterL\,  \TbcastI,\qquad m=1,\dots,r.
\end{align}

\subsection{Final cost model}
\label{sec:final}

The summary of considerations from Sections~\ref{sec:flops}--\ref{sec:comm} is to evaluate the
computational costs of the numerical method from Sections~\ref{sec:DGM}--\ref{sec:DDM} as the
sum of the wall-clock times for computation and communications over all time steps
$m=1,\dots,r$ as
\begin{align}
  \label{costsT}
  \costs := \sum\nolimits_{m=1}^r \costsm,\qquad
  \costsm := \left(\Wtimem + \Tcomm\right),\quad m=1,\dots,r,
\end{align}
where $\Tcomm :=\Tgathr + \Tbcast$, cf.~\eqref{addm8} and \eqref{comm4}.
The formula \eqref{costsT} is used in the following section to investigate the influence
of the number of subdomains and coarse elements on the wall-clock time of the computational
process. Moreover, it is employed in Section~\ref{sec:ADDM} for the prediction of the
computational costs for the next time steps after the re-meshing.

\subsection{Numerical study of computational costs}
\label{sec:num0}
We present the study of the efficiency of the two-level
hybrid Schwarz preconditioners by a benchmark from
\cite{GiraldoRostelli_JCP08,Robert_JAS93},
namely  the rising thermal bubble test case.
It simulates the evolution of a warm bubble in a constant potential temperature field.
The square computational domain is $\Om=(0,1000)\times(0,1000)\,\si{m}$, the final physical time
$T=50\,\si{s}$, and we prescribe 
the constant mean potential temperature $\bar{\tempP} = 300\,\si{K}$.
The Exner pressure follows from the hydrostatic balance \eqref{balance}
as $\pressE = 1 - g/(\ccP \bar{\tempP}) x_2$. The initial velocity is equal to zero.
Moreover, to the mean flow, we add 
the potential temperature perturbations
\begin{align}
  \label{thermal1}
  \tempP' = 
  \begin{cases}
    0 & \mbox{for } r > r_c, \\
    \tfrac12 \tempP_c \left(1 + \cos(\pi\, r / r_c)\right)    & \mbox{for } r \le r_c, \\
  \end{cases}
\end{align}
where $\tempP_c= 0.5\,\si{K}$, $r = | x - x_c|$, $x_c=(500, 300)\,\si{m}$ and $r_c = 250\,\si{m}$.
The no-flux boundary condition $\bkv\cdot\bkn=0$ is prescribed on the whole boundary.

For this numerical study we used a fixed (non-adapted)
quasi-uniform triangular grid with $\#\Thm =\numf{9074}$ elements
using $p=3$ polynomial approximation in space and $q=1$ in time. Therefore,
the resulting system $\Am$
has size $\Nhm=\numf{725920}$ and the number of non-zero entries $\NZm=\numf{230656000}$
for all $m=1,\dots,r$.
The time step is chosen adaptively according to \eqref{CT} with $c_T = 1$, which
leads to $r=1049$ time steps. The stopping condition for the
nonlinear algebraic solver \eqref{CA} was used with $c_A=10^{-3}$. Therefore,
the number of Newton and GMRES iterations is varying for $m=1,\dots,r$.


We carried out the computations using $\Mm=4,8,16,\dots,512$ subdomains and
the coarse mesh $\THm$ having $\Mm$, $2\Mm$, $4\Mm$, $\dots$, $32\Mm$ coarse elements
(i.e., with the number of splitting $\sm = 1,2,4,\dots,32$, cf.~Section~\ref{sec:model0}).
The number of elements of the local problems is 
$\#\Thmi \approx \#\Thm / \Mm$ for all $i=1,\dots,\Mm$.
Table~\ref{tab:fixed} presents the observed quantities in sum over
all time levels $m=1,\dots,r$.
Namely, we find the \\[1mm]
{\renewcommand{\arraystretch}{1.1}
\begin{tabular}{ll}
  $\iterNa=\sum\nolimits_{m=1}^r\iterN$ & Newton iterations,\\
  $\iterLa=\sum\nolimits_{m=1}^r\iterL$ & GMRES iterations, \\
  $\callCa=\sum\nolimits_{m=1}^r\callC$ & evaluations of matrix $\Am$,\\
\end{tabular}
\begin{tabular}{ll}
$\fl=\sum\nolimits_{m=1}^r\flm$ & floating point operations  \eqref{flopsF3}, \\
$\Wtime=\sum\nolimits_{m=1}^r\Wtimem$  & the wall-clock time of computation \eqref{addm8},\\
  $\comm=\sum\nolimits_{m=1}^r\commm$ &  the wall-clock time of
  communication \eqref{comm4}.\\
\end{tabular}
}
\\[1mm]
We emphasize that $\Wtime$ is the real time measured by
{\it MPI\_Wtime} function. On the other
hand, the communication $\comm$ represents the estimated wall-clock time for
$\gathr$ and $\bcast$
by \eqref{comm2} and \eqref{comm4} with fitted
parameters $\alpha_g$,  $\beta_g$,  $\gamma_g$,  $\alpha_b$, $\beta_b$ and $\gamma_b$ for
data from Section~\ref{sec:comm} since we have not enough cores for
the maximal number of subdomains $\Mm=512$.

\begin{table}
  \caption{Computational costs for hybrid Schwarz preconditioner
    for varying number of subdomains $\Mm$, number of elements of the local problem
    $n_i\approx \#\Thm$, and coarse mesh elements $n_0=\#\THm$.}
\label{tab:fixed}
  \begin{center}
    \setlength{\tabcolsep}{2.5pt}
    \begin{minipage}[t]{0.49\textwidth}
      \input{Figs_fixed/costsALL_N004.tex}\\[1mm]
      \input{Figs_fixed/costsALL_N008.tex}\\[1mm]
      \input{Figs_fixed/costsALL_N016.tex}
    \end{minipage}
    \hspace{0.01\textwidth}
    \begin{minipage}[t]{0.49\textwidth}
      \input{Figs_fixed/costsALL_N032.tex}\\[1mm]
      \input{Figs_fixed/costsALL_N064.tex}\\[1mm]
      \input{Figs_fixed/costsALL_N128.tex}\\[1mm]
      \input{Figs_fixed/costsALL_N256.tex}\\[1mm]
      \input{Figs_fixed/costsALL_N512.tex}
    \end{minipage}

  \end{center}
\end{table}




We observe that increasing $\#\THm$ reduces the
number of GMRES iterations $\iterLa$ since 
the finer coarse mesh provides more information.
On the other hand, the number of $\fl$ is decreasing for
increasing $\#\THm$ only if the number of coarse elements $\#\THm$ is
at least several times smaller than the number of elements of the local problems $\#\Thmi$.
Otherwise, the computational costs of the solution of the coarse problem are non-negligible
and they prolong the computation.

Moreover, we find that the minimum in terms of $\fl$ ($\Mm=\#\THm = 64$) does not correspond to the
minimum in terms of $\Wtime$ ($\Mm=\#\THm = 128$) since the speed of the computations
depends also on the size of the system, see Figure~\ref{fig:speed0}. We note that these
results are hardware/software dependent.
Obviously, the wall-clock time  of communication $\comm$ is increasing for increasing
$\Mm$ since there are more cores which communicate, and similarly,
it is increasing for the size of the coarse problem since the transmitted vectors are longer.
Finally, we conclude that the wall-clock time
for communication is several times smaller than the time for computation.
Nevertheless, it is not negligible.

\section{Adaptive domain decomposition}
\label{sec:ADDM}

The goal of this section is to present an adaptive domain decomposition technique,
which (approximately) solves Problem~\ref{prob:main}.
For the generated sequence of systems \eqref{addm1} in
time layers $m=1,\dots,r$, the aim is
to choose the number of subdomains $\Mm$ and splitting $\sm$ 
such that the computational costs are minimal.
We recall that we consider four contributions to the computational costs:
the factorization and substitution of the local and global coarse systems \eqref{flops0}.

We employ the computational cost model \eqref{costsT} composed from \eqref{addm8}
and \eqref{comm4},
where we assume that $\fl$, $\speed$, and $\comm$ can be estimated from the empirical
formulas \eqref{addm5}, \eqref{addm7}, and \eqref{comm2},  respectively.
The parameters $c$, $C$, $\mu$, $\nu$, $a$,  $b$ and $\Smax$ from these formulas
are evaluated by the standard least square fitting employing the available information
from previous time steps which have to be stored during the time marching computation.
However, the local problems can be factorized and assembled in parallel.
Hence, it is sufficient to store the information only about the largest one, we index it by $i=1$.
The information about the coarse problem has index $i=0$ in the following.

Let $m$ be the current time step and $\mlev\ge2$.
For the last $\mlev$ adaptation cycles, we store the following
quantities for all $k=1,\dots,\mlev$ (in agreement with previous notation)
\begin{itemize}
  \itemsep0pt 
\item $\Nhmk$ -- the size of the global algebraic system \eqref{addm1},
\item $\Nhhmk$ -- the size of the largest local fine system \eqref{flops0},
i.e., $\Nhhmk =\max_{i=1,\dots,\Mmk} \dim (\bSmki)$,
\item $\NHmk$ -- the size of the global coarse system \eqref{flops0},
\item $\flfacmi$, $\flassmi$ -- the corresponding values of {\flops}
  for factorization and substitution of local ($i=1$) and global coarse ($i=0$) problems
  given by MUMPS, cf.~Section~\ref{sec:flops},
\item $\spfacmi$, $\spassmi$ -- the average {\it speed}
  for factorization and substitution of local ($i=1$) and global coarse ($i=0$) problems
  in time step $m-k$,
 \item $\callCk$ -- the number of updates of matrix $\Amk$ in time step $m-k$,
 \item $\iterLk$ -- the number of GMRES iteration in time step $m-k$ (sum of all Newton steps),
 \item $\Mmk$ - the number of subdomains (= the number of used cores) in time step $m-k$,
 \item $\TgathrIk$, $\TbcastIk$ -- the average wall-clock time of one $\gathr$ and $\bcast$
   in time step $m-k$.
\end{itemize}
If a mesh $\Thm$ is used for several time steps, then we store the previous data
only for the last time step executed  on $\Thm$ and ignore the others. We do not
emphasise this selection in the notation explicitly.

\subsection{Adaptive domain decomposition algorithm}
\label{sec:ADDM_alg}
We are ready to introduce the desired adaptive algorithm.
As mentioned above, first we set the parameters of empirical formulas using the stored
data and them among all candidates of $\Mm$ and $\sm$, we choose those which minimize
the computational cost model. Particularly, the algorithm has the following steps:
\begin{enumerate}[label=(A{\arabic*}),ref=(A\arabic*)]
\item {\em fitting of parameters}
  \begin{enumerate}[label=(A\arabic{enumi}\alph*)] 
  \item {\flops} -- in virtue of \eqref{addm5}, we define the
    exponential functions estimating the number of {\flops} as
    \begin{align}
      \label{fit1}
      \Flfac^i(N)  &:= C_{i} \, N^{\mu_0}  ,\qquad
      \{C_i,\mu_i\} =\arg \min_{C,\mu} \sum\nolimits_{k=1}^\mlev \normP{\flfacmi - C (\Nhhmi)^\mu}{}{2},
      \qquad
      i=0,1\\
      \Flass^i(N)  &:= c_{i} \, N^{\nu_i}  ,\qquad
      \{c_i,\nu_i\} =\arg \min_{c,\nu} \sum\nolimits_{k=1}^\mlev \normP{\flassmi - c (\Nhhmi)^\nu}{}{2},
      \qquad i=0,1.
      \notag
    \end{align}
    \label{alg1a}
  \item {\em speed} of computation --
    in virtue of \eqref{addm7}, we define the functions estimating the speed
    of computations as 
    \begin{align}
      \label{fit2}
      \Spfac^i(N)  &:= \min(\Smax_{\fac},\ a_{\fac,i} \, N + b_{\fac,i}),\quad
      \{a_{\fac,i}, b_{\fac,i} \} =\arg \min_{a,b} \sum\nolimits_{k=1}^\mlev \normP{\spfacmi - a\, \Nhhmi - b}{}{2},
      \\
      \Spass^i(N)  &:= \min(\Smax_{\ass},\ a_{\ass,i} \, N + b_{\ass,i}),\quad
      \{a_{\ass,i}, b_{\ass,i} \} =\arg \min_{a,b} \sum\nolimits_{k=1}^\mlev \normP{\spassmi - a\, \Nhhmi - b}{}{2},
      \notag \\
      & \mbox{where} \quad
      \Smax_{\fac}=\max_{k=1,\dots,m-1,\ i=0,1} \spfacmi, \quad
      \Smax_{\ass}=\max_{k=1,\dots,m-1,\ i=0,1} \spassmi\qquad \mbox{for } i=0,1. \notag
    \end{align}
    \label{alg1b}
  \item  {\em number of updates and iterations} -- the predictions of the
    number of  updates of matrix $\Amk$  and
    the number of GMRES iterations in the Newton solver in the next time step $m$ is difficult
    to predict, so we use just an average of several last time steps, i.e.,
    \begin{align}
      \label{fit3}
      \bcallC := \frac{1}{\mlev} \sum\nolimits_{k=1}^\mlev \callCk,\qquad
      \biterL := \frac{1}{\mlev} \sum\nolimits_{k=1}^\mlev \iterLk,
    \end{align}
    \label{alg1c}
  \item {\em communication} times -- in virtue of \eqref{comm2}, we define the functions
    approximating the communication costs (including latency and bandwidth) by
    \begin{align}
      \label{fit4c}
      \TComg(N,P) & = \alpha_g \log_s{P} + \beta_g N + \gamma_g, \\
      & \{\alpha_g,\beta_g,\gamma_g\} = \arg \min_{\alpha,\beta,\gamma}
      \sum\nolimits_{k=1}^\mlev \normP{\alpha \log_2{\Mmk} + \beta \Nhmk + \gamma - \TgathrIk}{}{2},
      \notag \\
     \TComb(N,P) &= \alpha_b  \log_s{P}  + \beta_b N + \gamma_b, \notag \\
      & \{\alpha_b,\beta_b,\gamma_b\} = \arg \min_{\alpha,\beta,\gamma}
      \sum\nolimits_{k=1}^\mlev \normP{\alpha \log_2{\Mmk} + \beta \NHmk + \gamma - \TbcastIk}{}{2},
      \notag
    \end{align}
    \label{alg1d}
  \end{enumerate}
\item let $\Thm$ and $\bSm$ be the mesh and functional space generated by the anisotropic
  $hp$-mesh adaptation technique, the size of the corresponding algebraic system
  is $\Nhm=\dim\bSm$. For all admissible number of subdomains $\Mm=1,2,\dots$ and
  number of splittings $\sm=1,2,\dots$, (cf. Remark~\ref{rem:step_all})
  \label{alg:step_all}
  \begin{enumerate}[label=(A\arabic{enumi}\alph*)]
  \item   set the corresponding (approximate) sizes of local and global coarse systems
    \begin{align}
      \label{fit4}
      N^1 := \frac{\Nhm}{\Mm} \qquad \mbox{and} \qquad
      N^0 := \Nhm \,\frac{\#\THm}{\#\Thm} = \Nhm \,\frac{\sm\,\Mm}{\#\Thm},
    \end{align}
    respectively,
  \item in virtue of \eqref{addm8a}--\eqref{addm8} and \eqref{fit1}--\eqref{fit3},
    estimate the wall-clock time for the computation of the next time step as
    \begin{align}
      \label{fit5}
      \WWtime(\Mm, \sm) := \bcallC
      \max\left(\frac{ \Flfac^0(N^0) }{ \Spfac^0(N^0)}\ ,\ 
      \frac{ \Flfac^1(N^1) }{ \Spfac^1(N^1)} \right)
      + \biterL \left( \frac{\Flass^0(N^0) }{ \Spass^0(N^0)}
      + \frac{\Flass^0(N^0) }{ \Spass^0(N^0)} \right),
    \end{align}
    where $N^0$ and $N^1$ are given by \eqref{fit4}. We recall that the used hybrid method
    requires the solution of the global coarse system after the solution of the local ones,
    hence we have to sum both contributions in the second term in \eqref{fit5}.
  \item in virtue of \eqref{comm4} and \eqref{fit4c},
    estimate the wall-clock time for the communication of the next time step as
    \begin{align}
      \label{fit5a}
      \Comm(\Mm,\sm) := \biterL \left( \TComg(N^1,\Mm) +  \TComb(N^0,\Mm) \right).
    \end{align}
  \end{enumerate}
\item in virtue of \eqref{costsT}, \eqref{fit5}, and \eqref{fit5a},
  choose the optimal values $\MmO$ and $\smO$ that minimize the predicted computational costs
  as 
  \begin{align}
    \label{fit6}
    \left\{\MmO , \smO \right\} = \arg \min_{\Mm=1,2,\dots,\ \sm =1,2,\dots}
    \Big( \WWtime(\Mm, \sm) + \Comm(\Mm, \sm)\Big).
  \end{align}
\end{enumerate}

\begin{remark}
  \label{rem:step_all}
  In practical computations, the terms ``all admissible number of subdomains $\Mm=1,2,\dots$ and
  number of splittings $\sm=1,2,\dots$'' in step \ref{alg:step_all}, has natural restrictions.
  First, $\Mm$ is bounded by the number of available cores. In principle,
  it is possible to use higher number of subdomains $\Mm$ than the number of available cores
  but in this case the computational
  cost model has to be modified. Moreover, we can assume that each subdomain has
  at least several elements $n^{\min}$ (typically 10, 100, \dots, depending on the application).
  From this, we have the restriction
  \begin{align}
    \label{fit11}
    \Mm  \leq M_m^{\max} := \#\Thm\ / \ n^{\min}.
  \end{align}
  Moreover, based on previous considerations and numerical studies, it is natural
  to restrict to the case when the number of elements of the coarse mesh is lower than the number of
  elements of the local meshes, i.e., 
  \begin{align}
    \label{fit12}
    \#\THm = \sm\,\Mm  \leq \#\Thm\ / \ \Mm\qquad \Rightarrow \quad
    \sm  \leq s_m^{\max} := \#\Thm\ / \ \Mm^2.
  \end{align}
\end{remark}

We note that the algorithm described above cannot  be applied in the first $\mlev$ steps.
So we have to start the computational process in another way, typically we prescribe
the number of subdomains such that the size of a local problem is approximately
equal to the given value $\kk$, thus we put $\Mm \approx \Nhm / \kk $.

Finally, it is necessary to choose the number of time layers $\mlev$ used
for the fitting of parameters. The smaller value can better reflect the local in time
character of the problem considered but the fitting of parameters is less stable.
It is sensitive namely for the speed of computations and wall-time of communications.
Therefore, based on numerical
experiments, we use the value $\mlev=5$ for steps \ref{alg1a} and \ref{alg1c}.
On the other hand, we use all the available history for steps \ref{alg1b} and \ref{alg1d},
i.e., $\mlev=m$.

\begin{figure} [t]
  \begin{center}
    \includegraphics[width=0.49\textwidth]{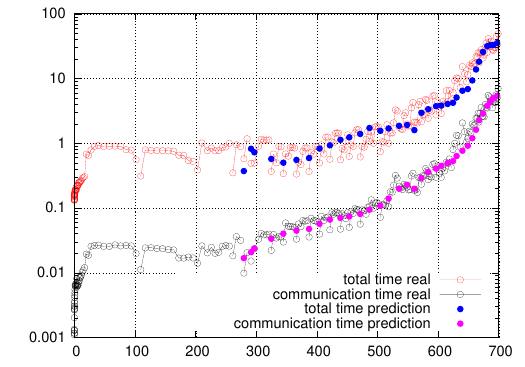 }  
    \includegraphics[width=0.49\textwidth]{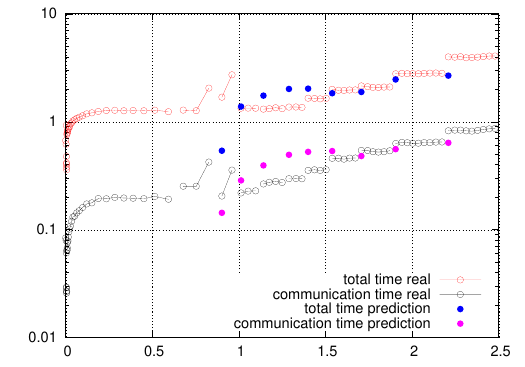}  
\end{center}
  \caption{Comparison of the predicted and real wall-clock times (vertical axes)
    for each time step (the horizontal axis is the physical time);
    examples from Section~\ref{sec:numerB} (left) and Section~\ref{sec:numerK} (right).
    } 
  \label{fig:predict}
\end{figure}

The quality of the prediction of the computational costs by technique \eqref{fit1}--\eqref{fit5a}
is shown in Figure~\ref{fig:predict} which compares the predicted and real wall-clock times
for the communications $\Tcomm$ and total time $\costsm$ (computations + communication)
for all time steps for examples from Sections~\ref{sec:numerB} and \ref{sec:numerK}.
Each node corresponds to one time step, the filled circles are the predicted times
for the next time step, and empty circles are the real times.  We observe a reasonably
good agreement
between the predicted and real times. We recall that the prediction of the computational costs
and the subsequent (``optimal'') domain decomposition is carried out only when a new mesh
is generated, and it can not be used in a few first time steps.

\section{Numerical experiments}
\label{sec:numer}

We present numerical examples, which demonstrate the efficiency
of the proposed adaptive domain decomposition algorithm from 
Section~\ref{sec:ADDM_alg}.
We always run the STDG technique with the anisotropic $hp$-mesh adaptation described
in Section~\ref{sec:STDG} with the same parameters, the linear algebraic systems
are solved iteratively by GMRES method preconditioned by the two-level hybrid Schwarz
preconditioner.
Note that their dimension grows during the time stepping.
We use three domain decomposition
variants of the choice $\Mm$ and $\sm$ at time level  $m=1,\dots,r$,
cf. Problem~\ref{prob:main}:
\begin{itemize}
\item [{\fix}] ({\em Fixed subdomains}):
  The number of subdomains $\Mm$ is fixed for all $m=1,\dots,r$,
  we set the values $\Mm=10$, 20, 30, 50, and 70.
  Additionally, $\Mm$ is chosen so that each $\Om_i$, $i=1,\dots,\Mm$
  has at least 10 elements, otherwise  $\Mm$ is decreased.
  The parameter $\sm$ is chosen such that $\#\THm \leq \tfrac12 \#\Thm / \Mm$, i.e, 
  the number of elements of the coarse grid is at most one half of the number of elements
  of each subdomain,  cf.~\eqref{fit12}.
\item [{\equi}] ({\em Fixed number of degrees of freedom}):
  The number of subdomains $\Mm$ is chosen such that each subdomain has approximately
  the same number of degrees of freedom (unknowns) $\kk$ for all time levels, i.e.,
  $\Mm = \Nhm/ \kk$ for all $m=1,\dots,r$, where $\Nhm=\dim\bSm$ is the size of the global system.
  We use the values $\kk=4000$, 6000, 8000, 10\,000, 12\,000, and 16\,000,
  and we require that $\Mm\ge2$. The parameter $\sm$ is chosen as in {\fix}.
\item [{\adapt}] ({\em adaptive choice})
  The number of subdomains $\Mm$ and the parameter $\sm$ are chosen adaptively
  using the algorithm from Section~\ref{sec:ADDM_alg}. Since it depends on the
  measuring of the wall-clock time during the computational process, which is a
  system function, we execute the computation with the same setting six times and denote these runs
  as {\exe} 1, \dots, {\exe} 6. For the first 5 steps, we prescribe $\Mm$ as in {\equi} with
  $\kk=2500$. 
\end{itemize}


In the following, we present the computations for two benchmarks
using the types of domain decomposition listed above and compare their computational costs.
Similarly  as in Section~\ref{sec:num0}, we investigate 
the total number of Newton iterations $\iterNa=\sum\nolimits_{m=1}^r\iterN$,
and total number of GMRES iterations $\iterLa=\sum\nolimits_{m=1}^r\iterL$. However, these
quantities have only informative character since the size of the systems $\Am$ is
changing for $m=1,\dots,r$.

\subsection{Rising thermal bubble}
\label{sec:numerB}

We consider the rising thermal bubble test case
\cite{GiraldoRostelli_JCP08,Robert_JAS93} which simulates the evolution of a warm bubble
in a constant potential temperature field. 
The setting is the same as in Section~\ref{sec:num0}, the only difference is that
the  final physical time is $T=700$\,\si{s}, and 
the mesh is not fixed but the anisotropic $hp$-mesh adaptation is used. 
The computational results for all three domain decomposition types 
{\fix}, {\equi}, and {\adapt} are given in Table~\ref{tab:adaptB}.
The columns $\Mm$ and $\#\THm$ denote the number of subdomains and the number of coarse
elements 
on the last time layer, respectively.
The next pair of columns shows the accumulated number of Newton ($\iterNa$), and
GMRES ($\iterLa$) iterations.
The next triple of columns shows the number of floating point operations ($\fl$),
wall-clock time of computation ($\Wtime$), and wall-clock time of communication ($\comm$),
the last column contains the total computational costs ($\costs= \Wtime + \comm$),
cf.~\eqref{costsT}.
The last three rows are the minimal, average, and maximal values of the {\adapt} technique.

\begin{table}
  \caption{Rising thermal bubble, comparison of computation costs for all domain
  decomposition variants {\fix}, {\equi}, and {\adapt}.}
  \label{tab:adaptB}
  \begin{center}
    \setlength{\tabcolsep}{4pt}
\input{Figs_adaptB/aDDM_costs.tex}

  \end{center}
\end{table}

We observe that the adaptive setting {\adapt} outperforms both other techniques.
For {\fix} and {\equi}, there exists a
setting 
giving costs which are not too far
from the costs of the {\adapt} technique. However, we do not know the optimal setting a priori,
it has to be found experimentally, and
it is typically case-dependent. On the other hand,  {\adapt} technique is fully automatic.
Additionally, we observe some differences in {\adapt} technique with
execution {\exe} 1, \dots, {\exe} 6. As explained above, the differences are caused by
measuring the wall-clock time during the computational process and therefore different
domain decomposition is carried out. However, the differences are in the order of a few percents.

\begin{figure} [t]
  \begin{center}
    \includegraphics[width=0.49\textwidth]{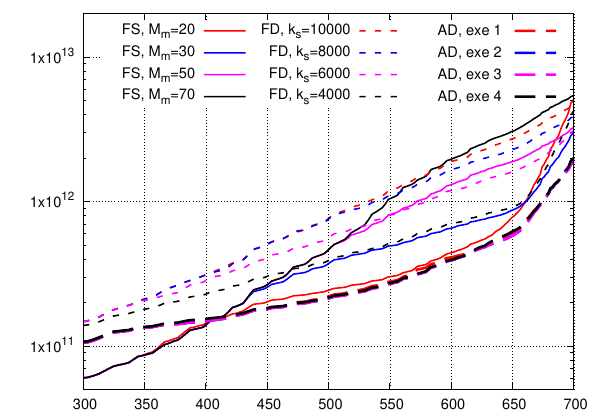}  
    \includegraphics[width=0.49\textwidth]{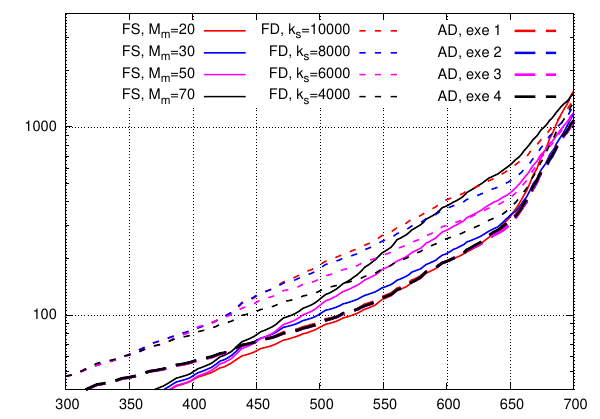}  
\end{center}
  \caption{Rising thermal bubble, comparison of all domain decomposition variants
    {\fix}, {\equi}, and {\adapt}, 
    the accumulated $\fl(s)$ (left) and total measured time $\costs(s)$ (right)
    for $t_s\in[300,T]$.} 
  \label{fig:adaptB1}
\end{figure}

\begin{figure} [t]
  \begin{center}
    \includegraphics[width=0.49\textwidth]{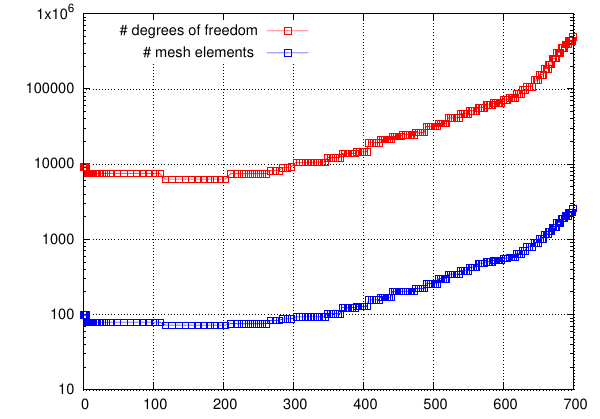}  
    \includegraphics[width=0.49\textwidth]{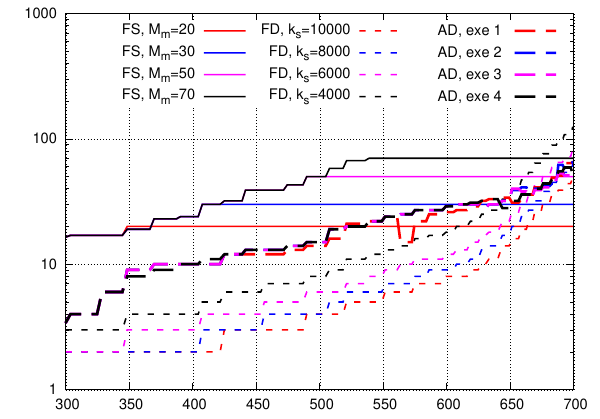}  
\end{center}
  \caption{Rising thermal bubble, 
    the number of elements $\#\Thm$ and degrees of freedom $\Nhm$ ($=\dim\bSm$)
    for each time step $t_m\in[0,T]$ (left),
    and the number of subdomains $\Mm$ for all variants  {\fix}, {\equi}, and {\adapt},
    $t_m\in[300,T]$ (right).}
  \label{fig:adaptB2}
\end{figure}

Moreover, Figure~\ref{fig:adaptB1} shows
dependence of the accumulated $\fl(s) = \sum_{m=1}^s\flm$ and
$\costs(s) = \sum_{m=1}^s\costsm$   with respect to physical time level $t_s\in[0,T]$,
cf.~\eqref{flopsF3} and \eqref{costsT}, for selected cases.
To better see the differences, we show only 
the time interval starting from $t=300$.
Furthermore, Figure~\ref{fig:adaptB2} shows the increase of
the number of elements $\#\Thm$, the number of degrees of freedom $\Nhm$ ($=\dim\bSm$) 
and the number of subdomains $\Mm$ for each time level $t_m\in[0,T]$.
The right-hand side graph in Figure~\ref{fig:adaptB2} illustrates the increase of
the number of subdomains $\Mm$. While it is given a priori for {\fix} and {\equi},
the {\adapt} technique is automatic, and it gives a slower increase in comparison
with the {\equi} choice. Moreover, small differences among executions {\exe} 1, \dots, {\exe} 4 are
observed.
We note that the number of elements and the number of degrees of freedom can
also decrease (see the left part of Figure~\ref{fig:adaptB2} for $t\leq 200\,\si{s}$) since
the mesh optimization process minimizes the number of degrees of freedom while keeping the accuracy
below the given tolerance, cf.~Section~\ref{sec:AMA} and the references therein.
For $t\leq 200\,\si{s}$, the shape of the bubble almost does not change.

\begin{figure} [t]
  \begin{center}
    \includegraphics[width=0.195\textwidth]{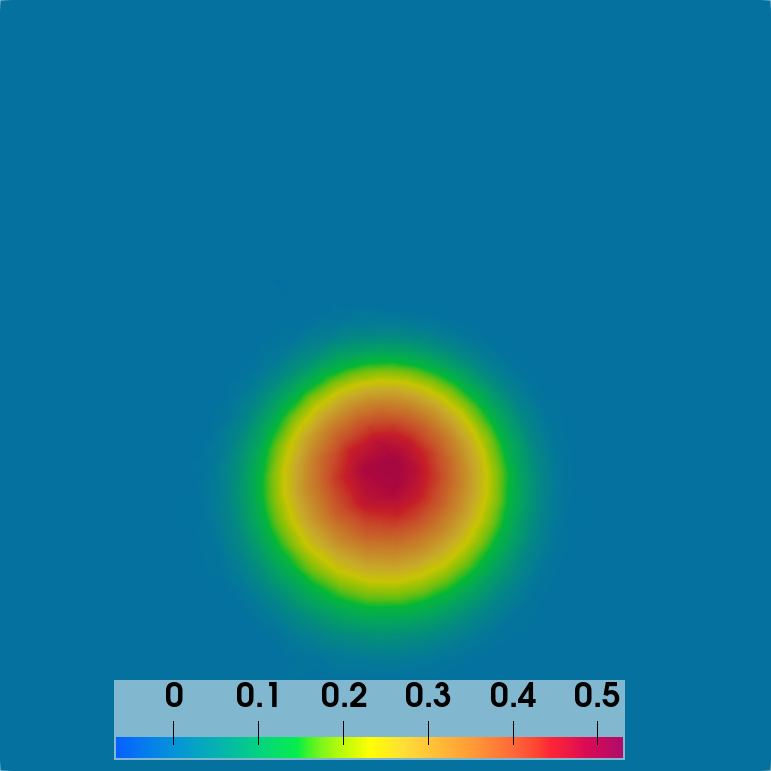}
    \includegraphics[width=0.195\textwidth]{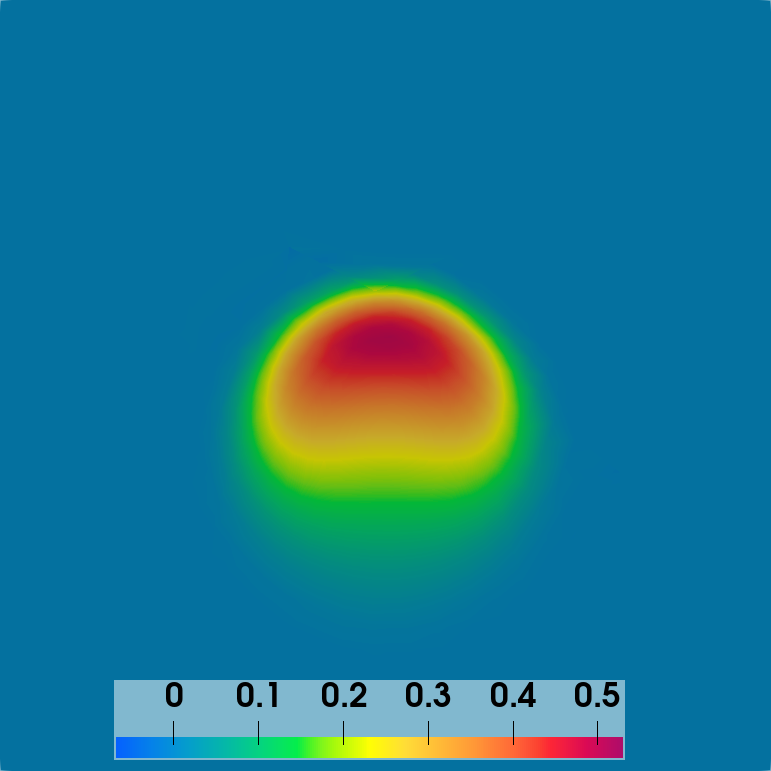}
    \includegraphics[width=0.195\textwidth]{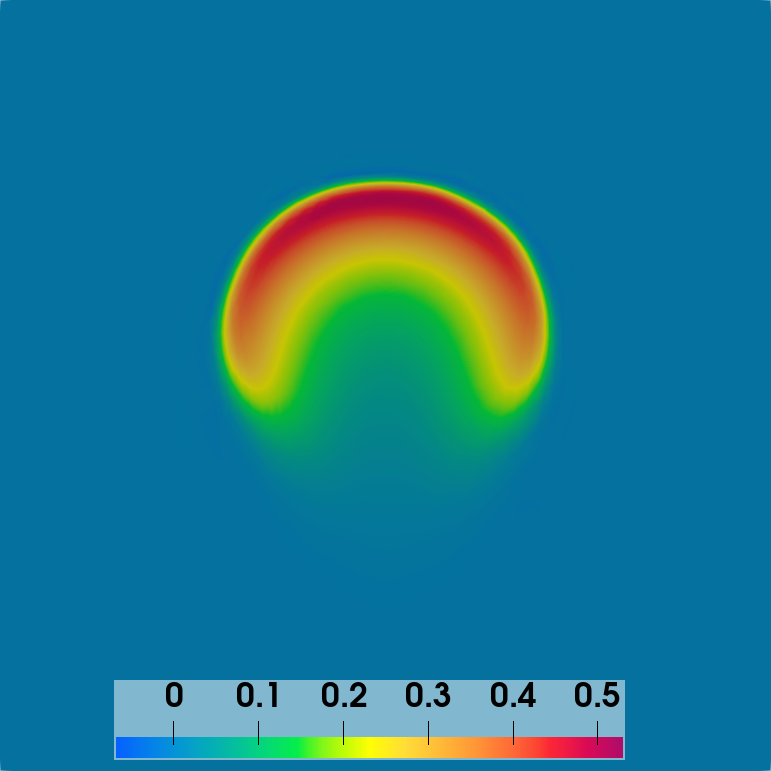}
    \includegraphics[width=0.195\textwidth]{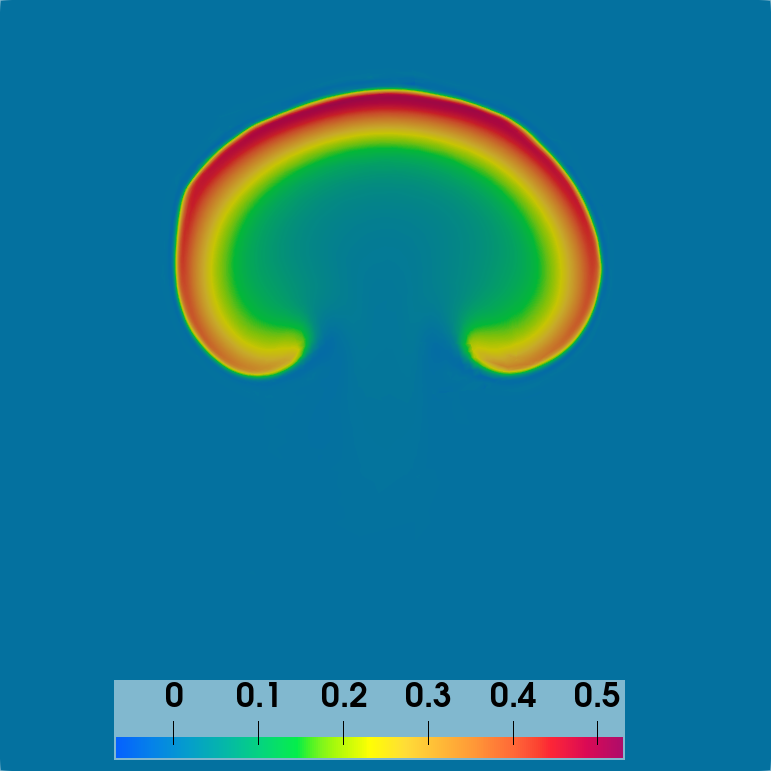}
    \includegraphics[width=0.195\textwidth]{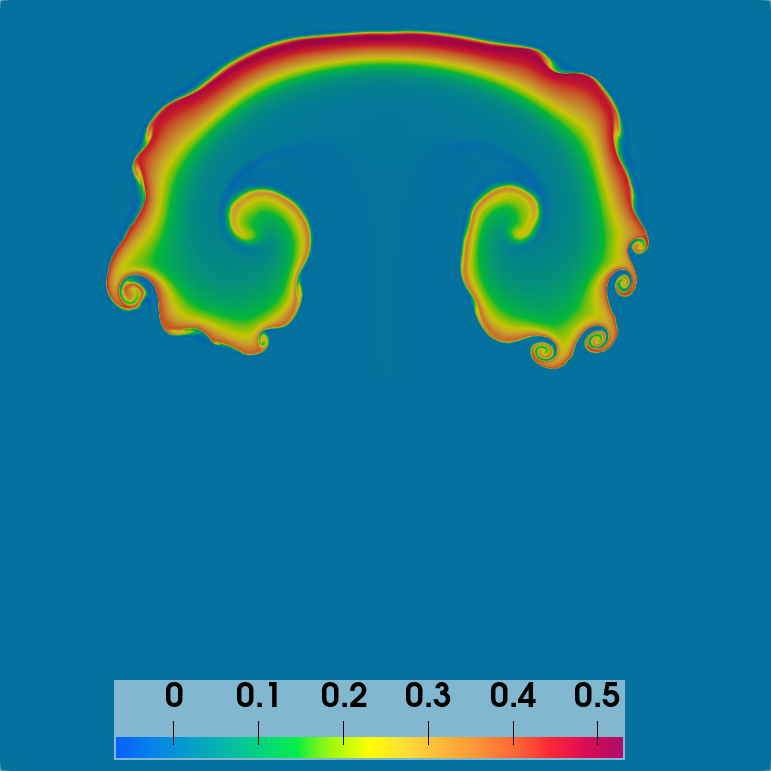}

    \includegraphics[width=0.195\textwidth]{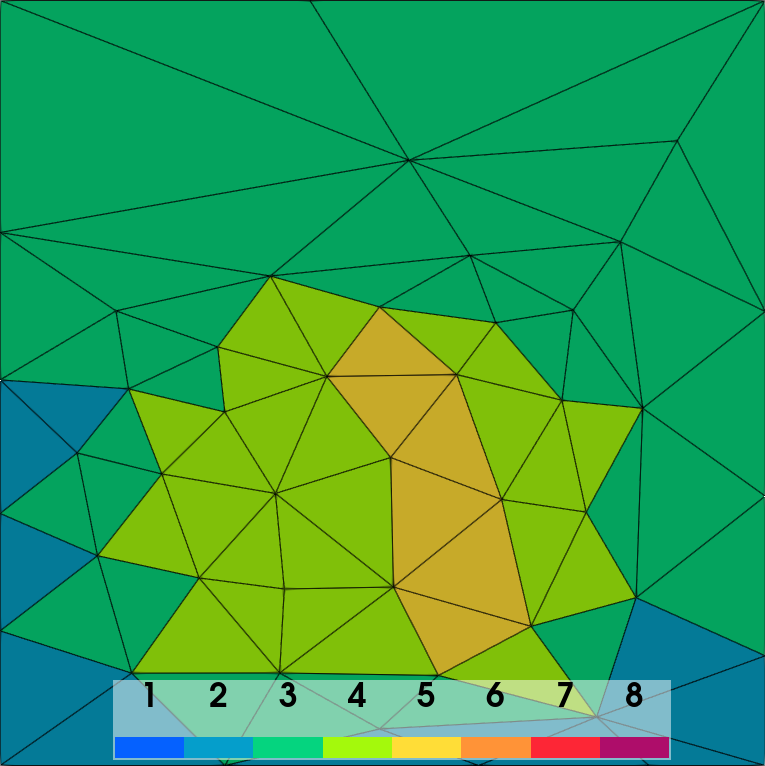}
    \includegraphics[width=0.195\textwidth]{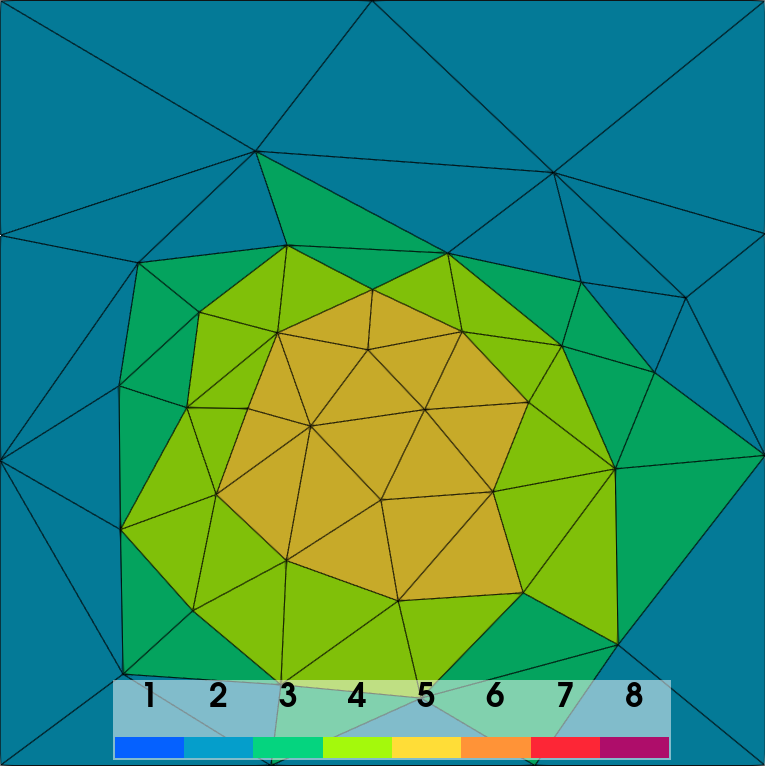}
    \includegraphics[width=0.195\textwidth]{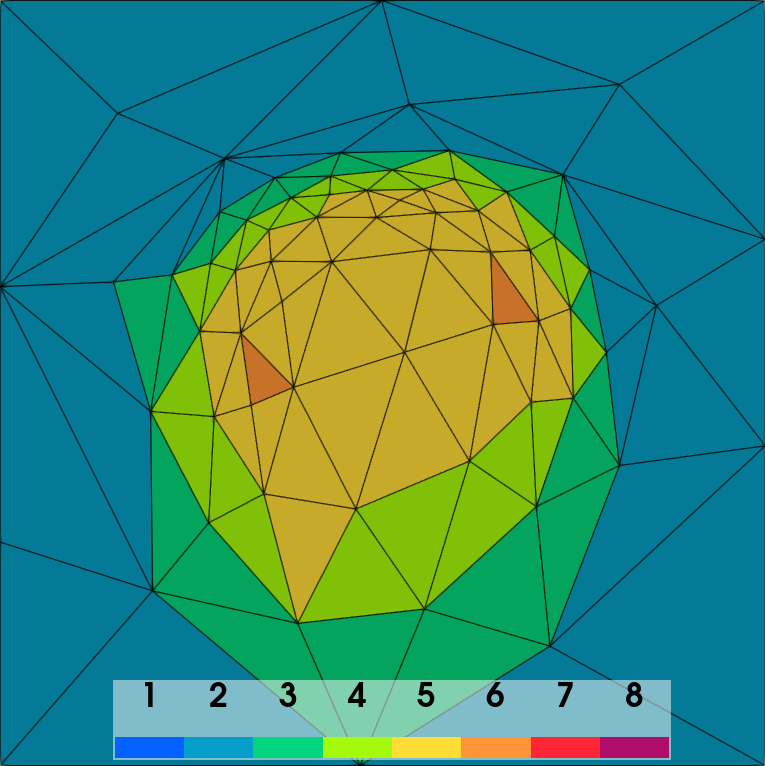}
    \includegraphics[width=0.195\textwidth]{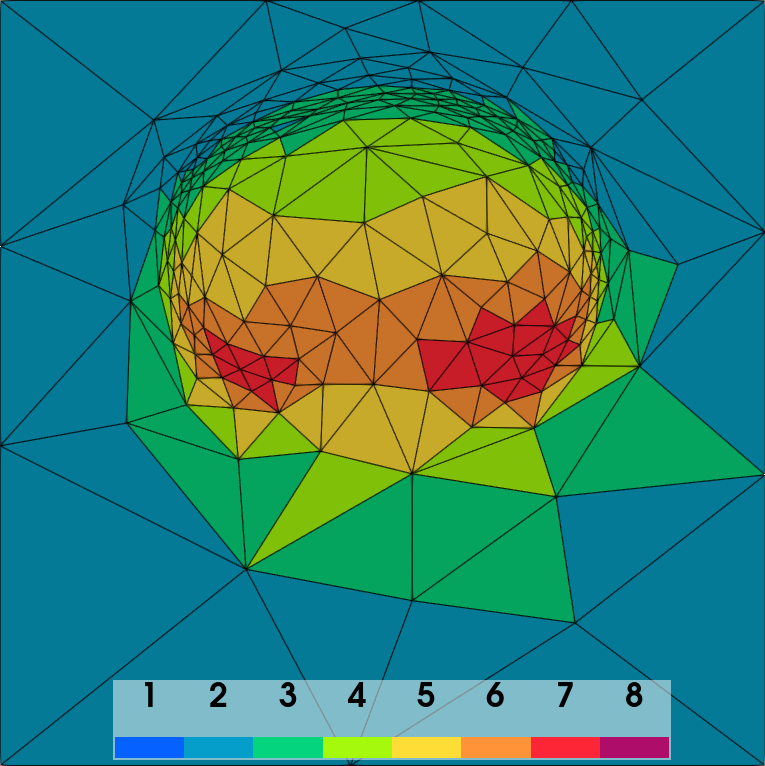}
    \includegraphics[width=0.195\textwidth]{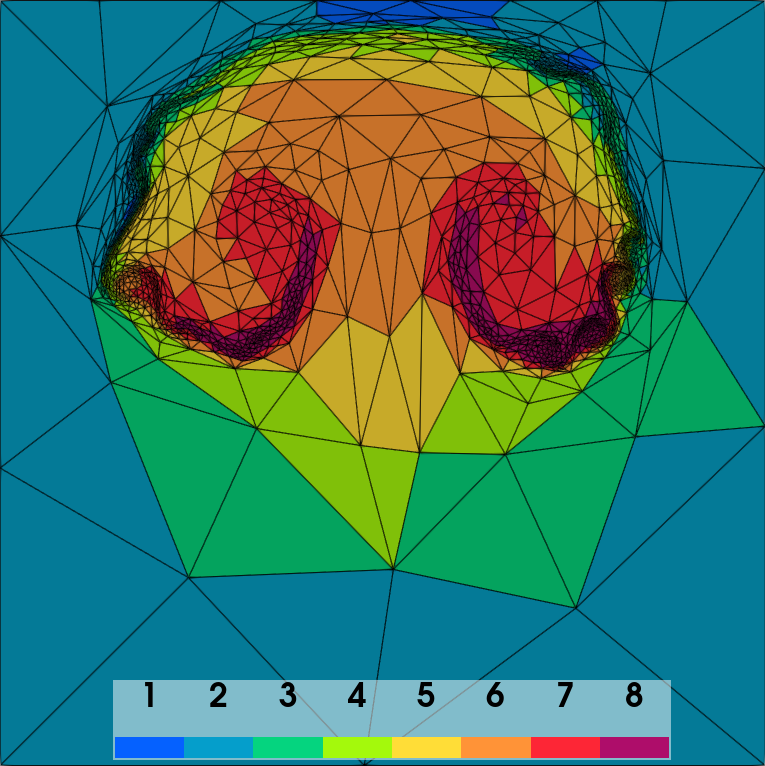}

    \includegraphics[width=0.195\textwidth]{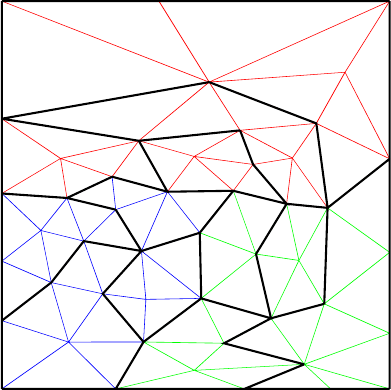}
    \includegraphics[width=0.195\textwidth]{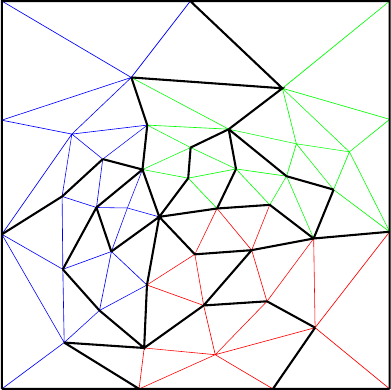}
    \includegraphics[width=0.195\textwidth]{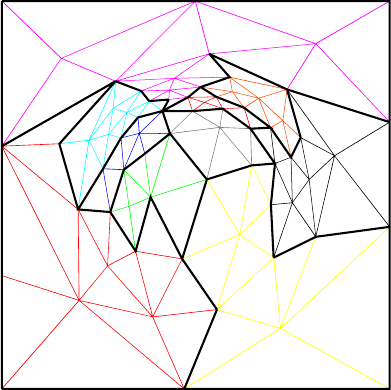}
    \includegraphics[width=0.195\textwidth]{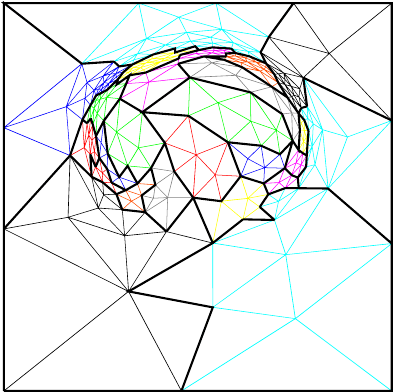}
    \includegraphics[width=0.195\textwidth]{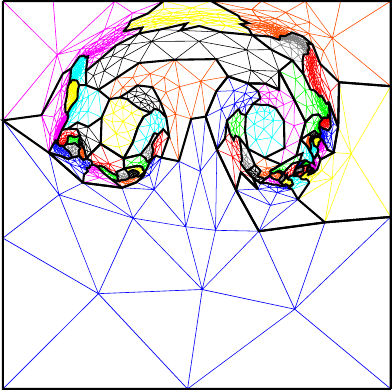}
\end{center}
  \caption{Rising thermal bubble, snapshots of the time-dependent simulation,
    the distribution of the potential temperature (first row), the anisotropic $hp$-mesh
    (second row) and the domain decomposition (third row) at time
    $t=$100, 250, 400, 550, and 700\si{s} (from left to right).}
  \label{fig:adaptB3}
\end{figure}

Finally, to illustrate the presented adaptive domain decomposition method,
Figure~\ref{fig:adaptB3} shows the snapshots of the time-dependent simulation
obtained by the {\adapt} technique. 
We present the distribution of the potential temperature $\tempP$ (cf.~\eqref{eq:tempP}),
the corresponding anisotropic $hp$-mesh  and the domain decomposition (the subdomains are colored)
at several time instants.
We observe the increasing flow complexity and the corresponding mesh refinement, and
increasing number of subdomains $\Mm$ for selected $m=1,\dots,r$.

\subsection{Kelvin-Helmholtz instability}
\label{sec:numerK}

The next example deals with  the simulation of the Kelvin-Helmholtz instability,
which appears when a velocity difference across the interface
between two fluids is presented, cf.~\cite{Svard2018,SCHROEDER2019} for example.
Then typical eddies appear, and they are challenging to capture numerically.
The problem is described by \eqref{eq:NS} with no gravity terms.

The computational domain
$\Om=(0,3)\times(0,1)\si{m}$ is extended
periodically in both directions and the final time is $T=2.5\si{s}$.
The initial conditions are given by
\begin{align}
  \label{kelvin1}
  \begin{cases}
    \rho = 2,\ v_1 = - 0.5 + \ve,\ v_2 =  \ve,\ \press = 2.5 &
    \mbox{ if } 0.25 < x_2 < 0.75, \\
    \rho = 1,\ v_1 =   0.5 + \ve,\quad  v_2 =  \ve,\ \press = 2.5 &
    \mbox{ if } 0.25 \ge x_2 \mbox{ or } x_2 \ge 0.75, \\
  \end{cases}
\end{align}
where we set $\ve= 0.01\sin(4\pi x_1)$ to trigger the instability. We consider the fluid viscosity
$\tilde{\mu} = 2\cdot10^{-4}\,\si{kg\cdot m^{-1}\cdot s^{-1}}$,
the heat capacity at constant pressure $c_p = 1005\,\si{J\cdot kg^{-1}\cdot K^{-1}}$,
the adiabatic Poisson constant $\gamma=1.4$, and the Prandtl number $Pr = 0.72$.

We compare all three domain decomposition types ({\fix}, {\equi}, and {\adapt}) 
in the same manner as in Section~\ref{sec:numerB}, the results 
are given in Table~\ref{tab:adaptK}.
We observe that {\adapt} is better than both other techniques
in terms of $\fl$, $\Wtime$, and $\costs$. The differences between {\adapt}  and the optimal
choices of {{\fix}} and {{\equi}} 
are not large but, once again, these optimal values have to be found experimentally
and differ from the optimal values for examples from Section~\ref{sec:numerB}.
The differences in {\adapt} technique with
execution {\exe}~1, \dots, {\exe}~6 in terms of $\Wtime$ are again about a few percents.

\begin{table}
  \caption{Kelvin-Helmholtz instability, comparison of computation costs for all domain
  decomposition variants {\fix}, {\equi}, and {\adapt}.} 
  \label{tab:adaptK}
  \begin{center}
    \setlength{\tabcolsep}{4pt}
\input{Figs_adaptK/aDDM_Kelvin_costs.tex}

  \end{center}
\end{table}

\begin{figure} [t]
  \begin{center}
    \includegraphics[width=0.49\textwidth]{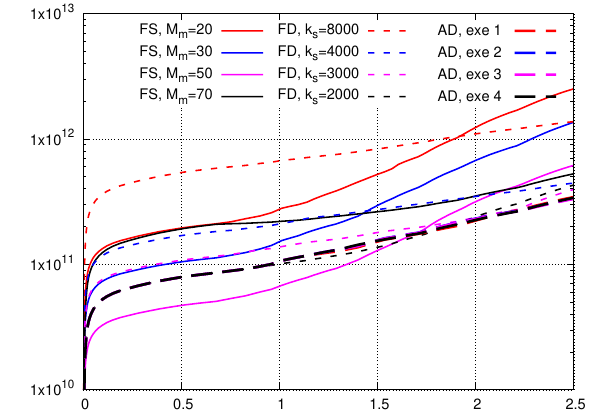}  
    \includegraphics[width=0.49\textwidth]{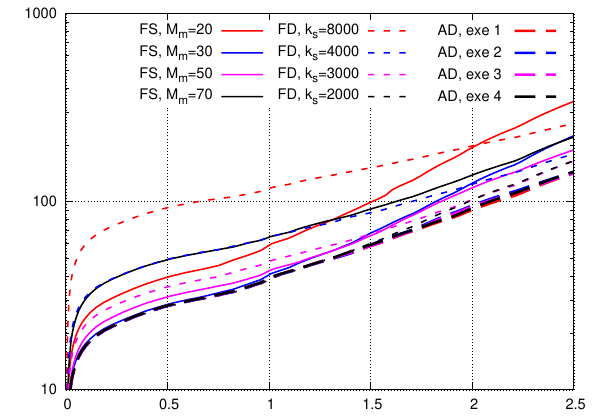}  

\end{center}
  \caption{Kelvin-Helmholtz instability,  comparison of all domain decomposition variants
    {\fix}, {\equi}, and {\adapt}, 
    the accumulated $\fl(s)$ (left) and total measured time $\costs(s)$ (right)
    with respect to physical time $t_s\in[0,T]$.}
  \label{fig:adaptK1}
\end{figure}

\begin{figure} [t]
  \begin{center}
    \includegraphics[width=0.49\textwidth]{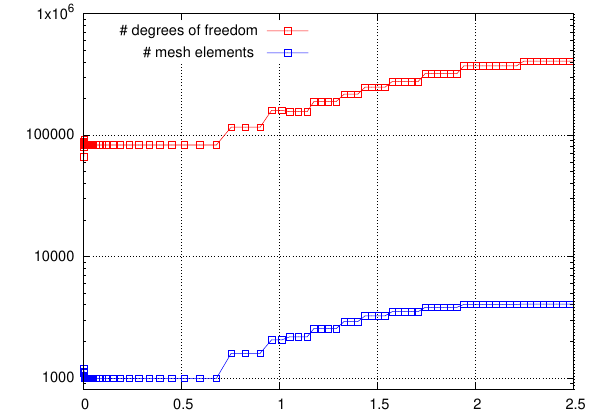}  
    \includegraphics[width=0.49\textwidth]{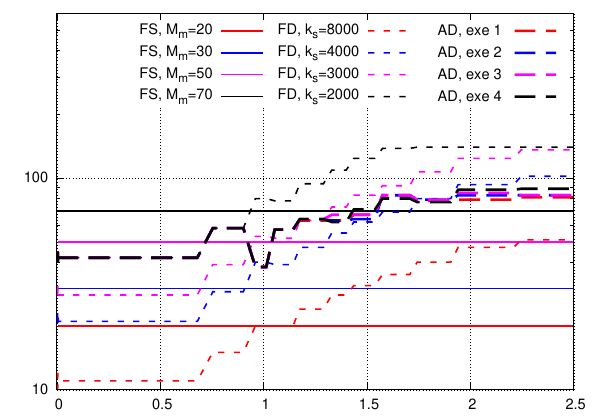}  

\end{center}
  \caption{Kelvin-Helmholtz instability,
    the number of elements $\#\Thm$ and degree of freedom $\Nhm$ ($=\dim\bSm$)
    for each time step $t_m\in[0,T]$ (left),
    and the number of subdomains $\Mm$  for all variants  {\fix}, {\equi}, and {\adapt},
    $t_m\in[0,T]$ (right).}
  \label{fig:adaptK2}
\end{figure}

Figure~\ref{fig:adaptK1} shows the dependence of the accumulated $\fl(s) = \sum_{m=1}^s\flm$ and
$\costs(s) = \sum_{m=1}^s\costsm$   with respect to physical time level $t_s\in[0,T]$,
cf.~\eqref{flopsF3} and \eqref{costsT}, for selected cases.
Furthermore, Figure~\ref{fig:adaptK2} shows the increase of
the number of elements $\#\Thm$, the number of degrees of freedom $\Nhm$ ($=\dim\bSm$) 
and the number of subdomains $\Mm$ for each time level $t_m\in[0,T]$.
The right-hand side graph in Figure~\ref{fig:adaptK2} illustrates the increase of
the number of subdomains $\Mm$.

\begin{figure} [t]
  \begin{center}
    \includegraphics[height=0.110\textwidth]{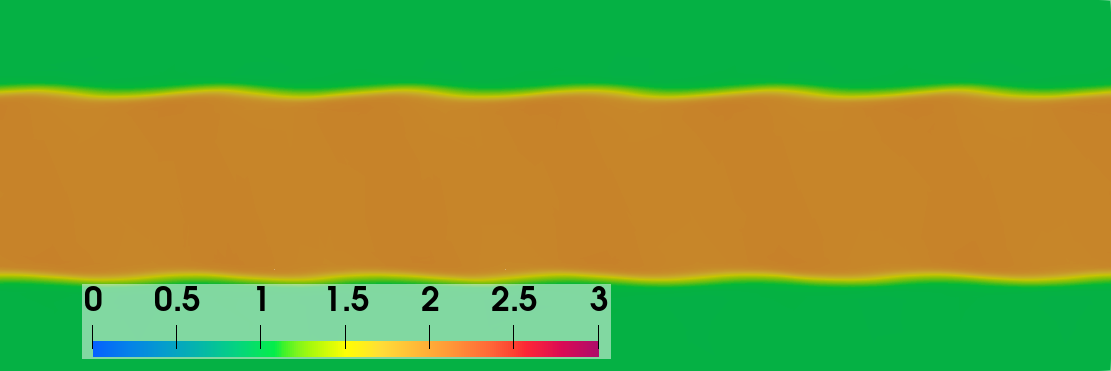}
    \includegraphics[height=0.110\textwidth]{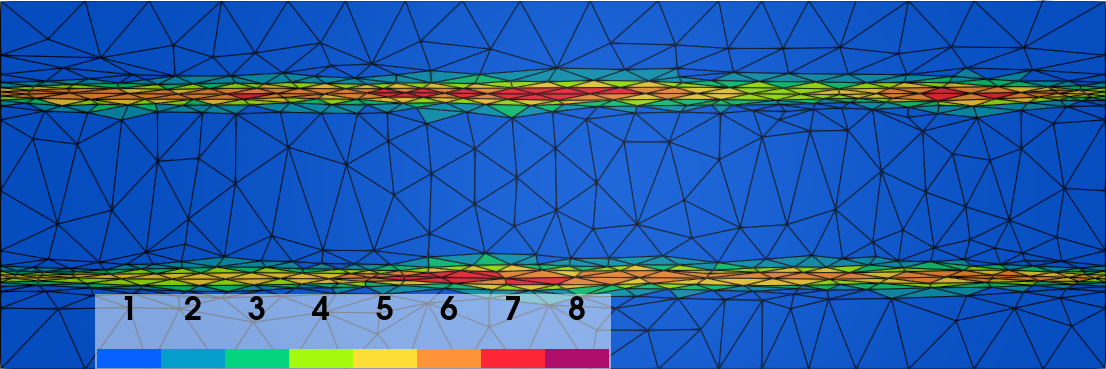}
    \includegraphics[height=0.110\textwidth]{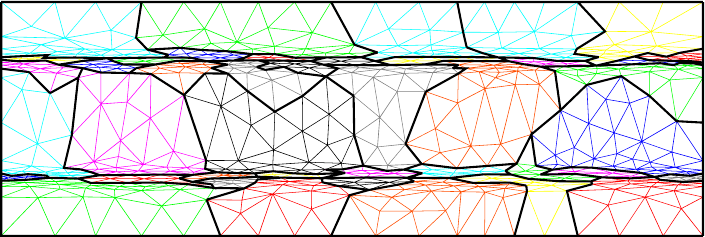}

    \includegraphics[height=0.110\textwidth]{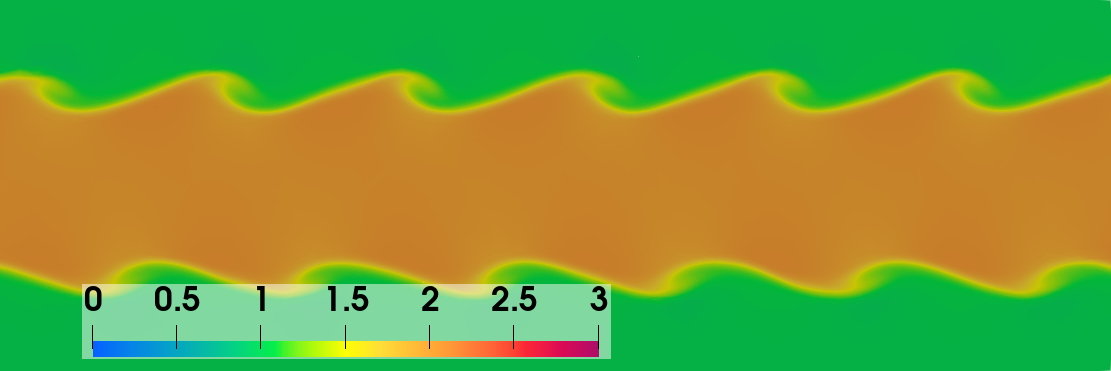}
    \includegraphics[height=0.110\textwidth]{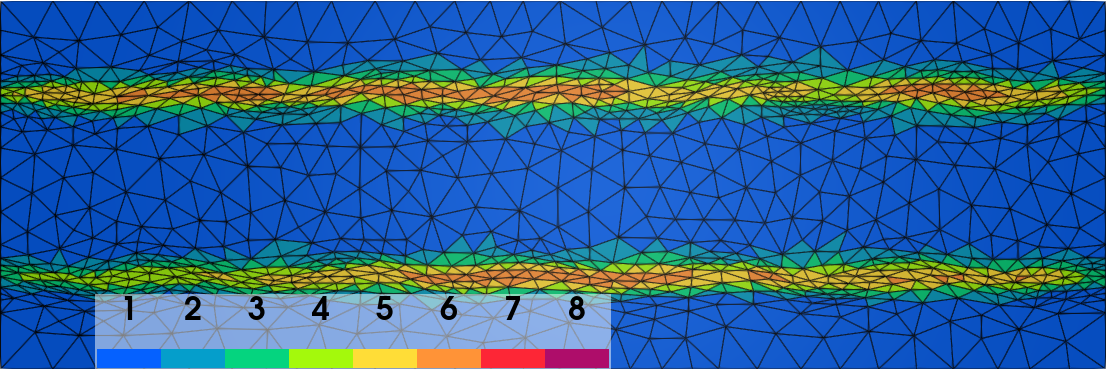}
    \includegraphics[height=0.110\textwidth]{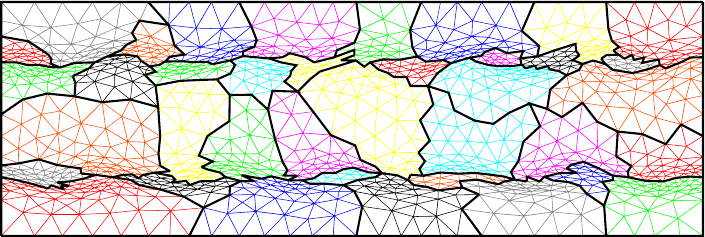}
    
    \includegraphics[height=0.110\textwidth]{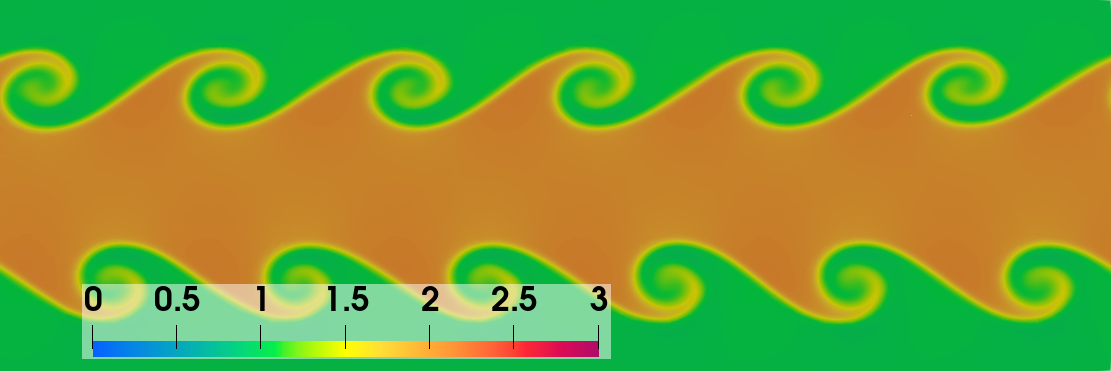}
    \includegraphics[height=0.110\textwidth]{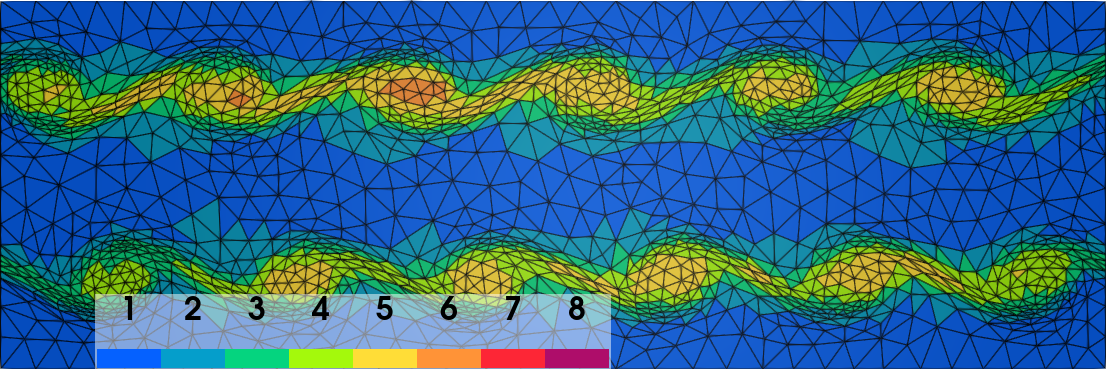}
    \includegraphics[height=0.110\textwidth]{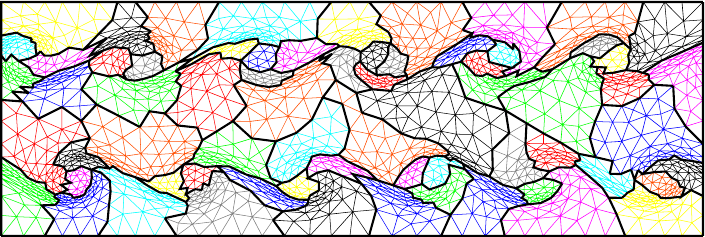}
    
    \includegraphics[height=0.110\textwidth]{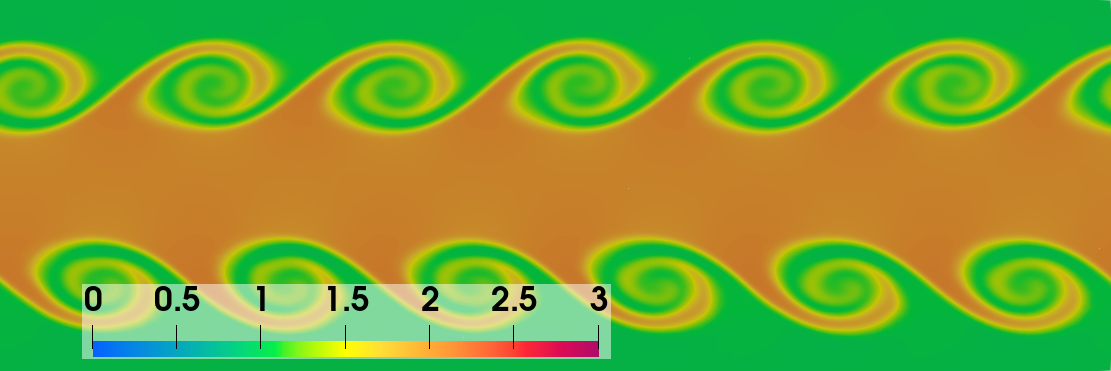}
    \includegraphics[height=0.110\textwidth]{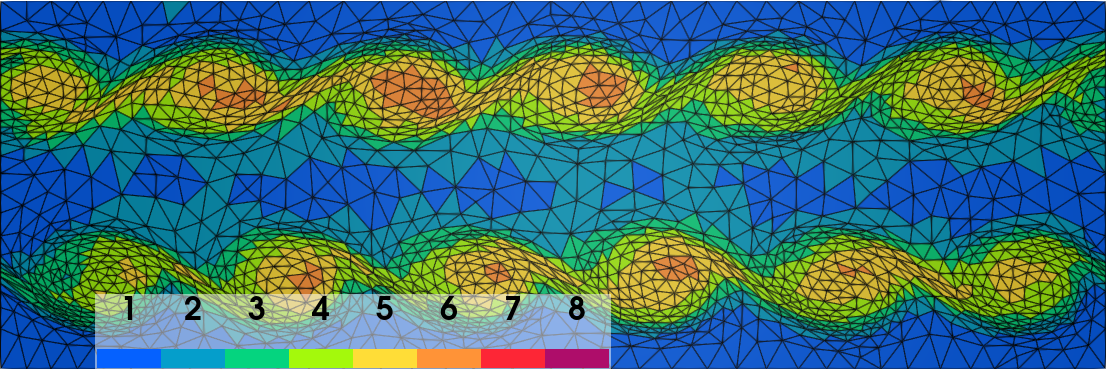}
    \includegraphics[height=0.110\textwidth]{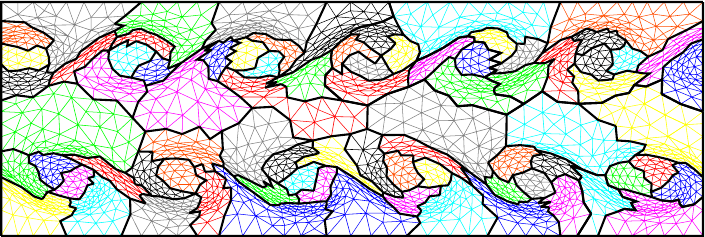}
    
    \includegraphics[height=0.110\textwidth]{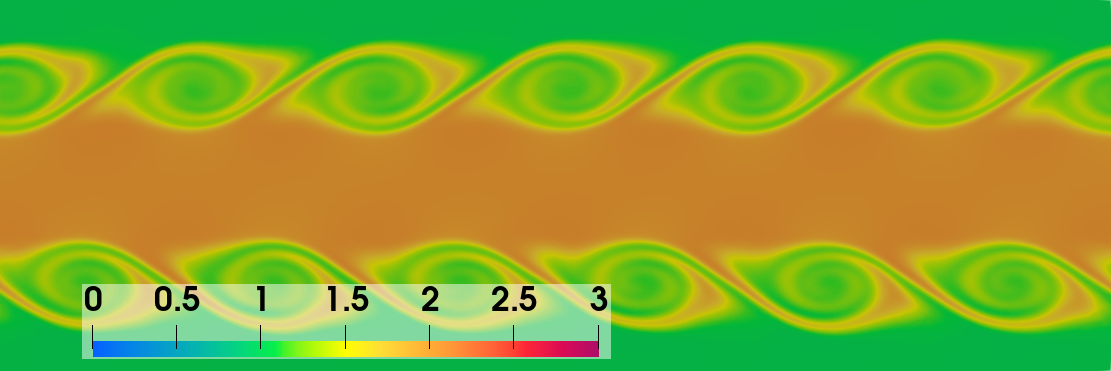}
    \includegraphics[height=0.110\textwidth]{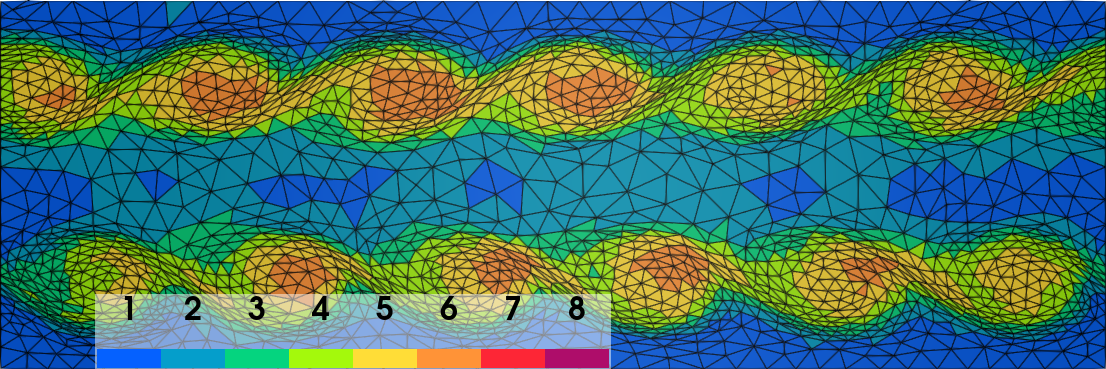}
    \includegraphics[height=0.110\textwidth]{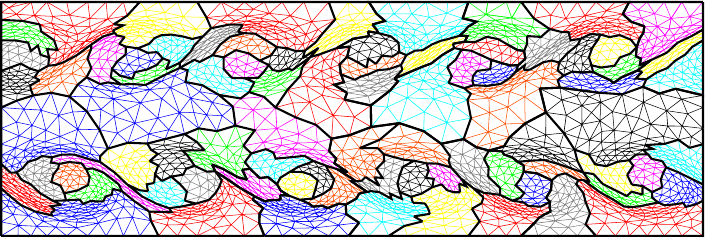}
  \end{center}
  \caption{Kelvin-Helmholtz instability,  snapshots of the time-dependent simulation,
    the distribution of the potential temperature (first column), the anisotropic $hp$-mesh
    (second column) and the domain decomposition (third column) at time
    $t=$0.5, 1.0, 1.5, 2.0, and 2.5\si{s} (from top to bottom).}
  \label{fig:adaptK3}
\end{figure}

Finally, Figure~\ref{fig:adaptK3} illustrates the presented adaptive domain decomposition method
obtained by the {\adapt} technique. 
We show the distribution of density,
the corresponding anisotropic $hp$-mesh, and the domain decomposition
at several time instants.
We observe the increasing flow complexity, the corresponding mesh refinement, and
the increasing number of subdomains $\Mm$ for selected $m=1,\dots,r$.
Remarkably, we observe several symmetries of the solution although the meshes
are not strictly symmetric.

\section{Conclusion}
\label{sec:concl}

We have dealt with the adaptive domain decomposition preconditioners for the solution
of algebraic systems arising from the space-time discontinuous Galerkin discretization.
In particular, we have considered two-level additive and hybrid Schwarz techniques, and proposed
a computational cost model involving the wall-clock time of computation and communications
of the preconditioner.
This model is at the core of an adaptive algorithm which selects the number of subdomains and
the  number of elements of the coarse mesh in order to  minimize the predicted
computational costs.
The efficiency of the approach was demonstrated on two well-known benchmarks
where the adaptive algorithm equipped with the cost model performed better than
the optimal versions
of the strategies with constant subdomain size and constant number of subdomains, setting of which
is even unknown beforehand.

The presented computational cost model can be extended to other types of discretization and
preconditioners in a straightforward way. Additionally, the model can be generalised to
include other parts
of the computational process. This will be the subject of our further research.

\subsection*{Data availability}

Data will be made available on request.

\subsection*{Conflict of interest}
The authors declare that they have no conflicts of interest.

\subsection*{Acknowledgment}
J.~Šístek acknowledges the support by grant No.~23-06159S of the Czech Science Foundation and
by the Institute of Mathematics of the Czech Academy of Sciences (RVO:67985840).
Both authors acknowledge also the membership in the Ne{\v c}as Center for Mathematical Modeling
ncmm.karlin.mff.cuni.cz.

\subsection*{Declaration of generative AI and AI-assisted technologies in the manuscript preparation process.}

During the preparation of this work the authors did not use any AI and AI-assisted technologies.


\end{document}

%% file: def.tex
\newtheorem{theorem}{Theorem}[section]

\theoremstyle{remark}
\newtheorem{remark}[theorem]{Remark}
\newtheorem{definition}[theorem]{Definition}
\newtheorem{problem}[theorem]{Problem}

\newcommand\bifont {\bm}

\newcommand\krt{{\mathscr T}}

\newcommand\bkn{\bifont{n}}
\newcommand\bku{\bifont{u}}
\newcommand\bkv{\bifont{v}}
\newcommand\bkx{\bifont{x}}
\newcommand\bky{\bifont{y}}
\newcommand\bkz{\bifont{z}}

\newcommand\cK {\mathcal{K}}

\newcommand\pd {\partial}
\newcommand{\pdt}[1]{ \pd_t{#1} }
\newcommand{\pdd}[2]{ \pd_{#2}{#1} }

\newcommand\R  {\mathbb{R}}
\newcommand\mN  {\mathbb{N}}

\newcommand\mI  {\pmb{I}}
\newcommand\ve {\varepsilon}
\newcommand\Om {{\Omega}}
\newcommand\oOm  {\overline{\Om}}
\newcommand\oO {\overline{\om}}

\newcommand\gom{{\Gamma}}

\newcommand\T{^{\mkern-1.5mu\mathsf{T}}}

\newcommand\dd {\mathrm{d}}
\newcommand\dt {{\,\dd t}}

\newcommand\Nhm {N_{m}}
\newcommand\Nhmi {{N_{m}^i}}
\newcommand\NHm {{N_{m}^0}}

\newcommand\Mm {M_{m}}
\newcommand\Mmk {{M_{m-k}}}
\newcommand\MmO {{M_{m}^{\mathrm{opt}}}}

\newcommand\Nhmk {N_{m-k}}
\newcommand\NHmk {N_{m-k}^0}
\newcommand\Nhhmk {N_{m-k}^1}
\newcommand\Nhhmi {N_{m-k}^i}

\newcommand\bbw {\bar{{\w}}}
\newcommand\wht {{\w}_{h\tau}}

\newcommand\whtm {{\w}_{h\tau}^m}
\newcommand\whtmM {{\w}_{h\tau}^{m-1}}

\newcommand\WmM {\bm{\xi}^{m-1}}
\newcommand\Wm {\bm{\xi}^{m}}
\newcommand\Wml {\bm{\xi}^{m}_\ell}
\newcommand\WmlM {\bm{\xi}^{m}_{\ell-1}}
\newcommand\WmN {\bm{\xi}^{m}_0}
\newcommand\dl {\bm{d}_\ell}

\newcommand\bW {\bar{\bm{\xi}}}

\newcommand\vv {\bifont{v}}

\newcommand\TT {\mathrm{T}}
\newcommand\Th {\krt_h}
\newcommand\Thm{{\krt_{h,m}}}
\newcommand\Thmi{\krt_{h,m}^i}
\newcommand\THm{{\krt_{H,m}}}

\newcommand{\norm}[2]{\left\|#1\right\|_{#2} }
\newcommand{\normP}[3]{\left\|#1\right\|_{#2}^{#3} }

\newcommand\ahL {{a_h^{\scriptscriptstyle\mathrm{L}}}}

\newcommand\ahm {{a_{h,m}}}
\newcommand\ahmL {{a_{h,m}^{\scriptscriptstyle\mathrm{L}}}}
\newcommand\Ahm {{A_{h,m}}}
\newcommand\AhmL {{A_{h,m}^{\scriptscriptstyle\mathrm{L}}}}

\newcommand\om{{\Omega}}

\renewcommand\bkn{\vec{n}}

\newcommand\etaAm {\eta^{\scriptscriptstyle\mathrm{A}}_{m}}
\newcommand\etaSm {\eta^{\scriptscriptstyle\mathrm{S}}_{m}}
\newcommand\etaTm {\eta^{\scriptscriptstyle\mathrm{T}}_{m}}

\newcommand\tjump[1]{\{\!\{{#1}\}\!\}}

\newcommand{\Lsp}[2]{\left({#1},{#2}\right)_{\Omega}}
\newcommand{\LSP}[3]{\left({#1},{#2}\right)_{#3}}

\newcommand\CA{{C_{{A}}}}
\newcommand\CT{{C_{{T}}}}
\newcommand\CL{{C_{I}}}

\newcommand\DoF{\mathrm{DoF}}

\newcommand\bg {{G}}
\newcommand\bk {{\bifont{k}}}

\newcommand\w {w}
\newcommand\bF {f}
\newcommand\bR {R}

\newcommand\press {\mathrm{p}}
\newcommand\pressE {\mathrm{P}}
\newcommand\ener {e}
\newcommand\ccP{c_{\press}}
\newcommand\ccV{c_{\mathrm{v}}}
\newcommand\stress{\tau}

\newcommand\temp{\theta}
\newcommand\tempP{\Theta}

\newcommand\pK {p_K}
\newcommand\ddK {d_K}
\newcommand\ddcK {d_{\cK}}

\newcommand\qK{p_{\cK}}

\newcommand\mbP{\bifont{P}} 
 
\newcommand\mbR{\bifont{R}}

\newcommand\bSm {S_{\!m}^{\!h,\tau}}

\newcommand\bSmhh {S_{\!m}^{\!h+,\tau}}
\newcommand\bSmtt {S_{\!m}^{\!h,\tau+}}

\newcommand\bSmi {S_{\!m,i}^{\!h,\tau}}
\newcommand\bSmki {S_{\!m-k,i}^{\!h,\tau}}
\newcommand\bSHm {S_{\!m}^{\!H,\tau}}
\newcommand\bSmM {S_{\!m-1}^{\!h,\tau}}
\newcommand\bS {S^{\!h\tau}}

\newcommand{\bvas} {\varphi}
\newcommand{\bVas} {\phi}
\newcommand{\bpsi} {\psi}
\newcommand{\bPsi} {\psi}

\newcommand{\balpha} {\pmb{\alpha}}
\newcommand{\bbeta} {\pmb{\beta}}

\newcommand \hpAMA{{$hp$-AMA}}

\newcommand \Bm{{\mathcal{B}_{h,m}}}

\newcommand \BHm{{\mathcal{B}_{H,m}}}

\newcommand \BBm{{B_{h,m}}}
\newcommand \BBHm{{B_{H,m}}}

\newcommand{\Fm} {\pmb{f}_m}

\newcommand{\bA} {\pmb{A}}
\newcommand{\Am} {\bA_m}
\newcommand{\Amk} {\bA_{m-k}}
\newcommand{\Ami} {\Am^i}
\newcommand{\Amij} {\Am^{i,j}}
\newcommand{\Amii} {\Am^{i,i}}
\newcommand{\AmN} {\Am^0}

\newcommand{\oOmi} {\bar{\Om}_{m}^{i}}
\newcommand{\Omi} {\Om_{m}^{i}}
\newcommand{\Omj} {\Om_{m}^{j}}

\newcommand\flops {{\em flops}}

\newcommand\fac{\mathrm{fac}}
\newcommand\ass{\mathrm{sub}}

\newcommand\speed {\mathsf{speed}}
\newcommand\Wtime {\mathsf{time}}
\newcommand\Wtimem {\Wtime^m}
\newcommand\WWtime {\mathsf{Time}}

\newcommand\spfaci {\speed_{\fac}^{m,i}}
\newcommand\spassi {\speed_{\ass}^{m,i}}

\newcommand\Wtfaci {\Wtime_{\fac}^{m,i}}
\newcommand\Wtassi {\Wtime_{\ass}^{m,i}}
\newcommand\WtassN {\Wtime_{\ass}^{m,0}}
\newcommand\fl {\mathsf{flops}}
\newcommand\Tfl {\mathsf{wt}}
\newcommand\Cfl {\mathsf{ct}}

\newcommand\flm {\fl^m}

\newcommand\flfaci {\fl_{\fac}^{m,i}}
\newcommand\flassi {\fl_{\ass}^{m,i}}
\newcommand\flassN {\fl_{\ass}^{m,0}}
\newcommand\FFfac {\fl_{\fac}^m}
\newcommand\FFass {\fl_{\ass}^m}

\newcommand\FFl {\mathsf{Fl}}
\newcommand\Flfac {\FFl_{\fac}}
\newcommand\Flass {\FFl_{\ass}}

\newcommand\flfacmi {\fl_{\fac}^{m-k,i}}

\newcommand\flassmi {\fl_{\ass}^{m-k,i}}

\newcommand\SSp {\mathsf{Sp}}
\newcommand\Spfac {\SSp_{\fac}}
\newcommand\Spass {\SSp_{\ass}}

\newcommand\spfacmi {\speed_{\fac}^{m-k,i}}

\newcommand\spassmi {\speed_{\ass}^{m-k,i}}

\newcommand\task {\mathrm{T}}

\newcommand\iter {\mathsf{it}}
\newcommand\iterN {\iter_m^{\mathrm{n}}}
\newcommand\iterL {\iter_m^{\mathrm{\ell}}}
\newcommand\iterLk {\iter_{m-k}^{\mathrm{\ell}}}
\newcommand\iterNa {\iter^{\mathrm{n}}}
\newcommand\iterLa {\iter^{\mathrm{\ell}}}
\newcommand\biterL {\overline{\iter}_m^{\mathrm{\ell}}}

\newcommand\costs {\mathsf{costs}}
\newcommand\costsm {\mathsf{costs}^m}
\newcommand\comm {\mathsf{com}}
\newcommand\Comm {\mathsf{Com}}
\newcommand\commm {\comm^m}

\newcommand\callC {\mathsf{F}_m^{\mathrm{A}}}
\newcommand\callCk {\mathsf{F}_{m-k}^{\mathrm{A}}}
\newcommand\callCa {\mathsf{F}^{\mathrm{A}}}
\newcommand\bcallC {\overline{\mathsf{F}}_m^{\mathrm{A}}}
\newcommand\NZm {Z_{m}}

\newcommand\RR{\mbR}
\newcommand\mRi{\RR_{m,i}}
\newcommand\mRiT{\RR_{m,i}\T}

\newcommand\PP{\mbP}

\newcommand\mPPi{\PP_{m,i}}

\newcommand\mRN{\RR_0}
\newcommand\mRNT{\RR_0\T}

\newcommand\mPPN{\PP_0}

\newcommand\NN{\pmb{N}}

\newcommand\add{\mathrm{add},2}

\newcommand\mPadd{\PP_{\add}}

\newcommand\mNadd{\NN_{\add}}

\newcommand\hy{\mathrm{hyb}}

\newcommand\mNhy{\NN_{\hy}}

\newcommand\mPhy{\PP_{\hy}}

\newcommand\sm{s_m}
\newcommand\smO {{s_{m}^{\mathrm{opt}}}}
\newcommand\Smax {{\bar{S}}}
\newcommand\fix {{FS}}
\newcommand\equi {{FD}}
\newcommand\adapt {{AD}}
\newcommand\kk {{k_{\mathrm s}}}
\newcommand\exe {{exe}}
\newcommand\bcast {\mathit{bcast}}
\newcommand\gathr {\mathit{allgather}}
\newcommand\Bcast {{bcast}}
\newcommand\Gathr {{allgathr}}

\newcommand\TbcastI {{\comm_{\Bcast}^{m,1}}}
\newcommand\TgathrI {{\comm_{\Gathr}^{m,1}}}
\newcommand\TbcastIk {{\comm_{\Bcast}^{m-k,1}}}
\newcommand\TgathrIk {{\comm_{\Gathr}^{m-k,1}}}
\newcommand\Tbcast {{\comm_{\Bcast}^m}}
\newcommand\Tgathr {{\comm_{\Gathr}^m}}
\newcommand\Tcomm {{\comm^m}}

\newcommand\TComg {{\Comm_{\Gathr}}}
\newcommand\TComb {{\Comm_{\Bcast}}}

\newcommand\mlev{\vartheta}


%% file: Figs_flops/speed.tex
 \begin{tabular}{cccc}
 \hline
  operation & grid & $a $ & $b $ \\
 \hline
    factorization             &     triangular                &   5.5827E+10 &   6.4325E+04 \\
    factorization             &     polygonal                 &   4.7066E+10 &   2.1832E+04 \\
    substitution              &     triangular                &   1.1884E+09 &   6.7913E+03 \\
    substitution              &     polygonal                 &   1.1009E+09 &   3.7104E+03 \\
 \hline
 \end{tabular}

%% file: Figs_fixed/costsALL_N004.tex
 {\footnotesize
 \begin{tabular}{ccc|ccc|ccc}
 \hline
  $ \Mm $  & $n_i$ & $n_0$ 
  & $ \iterNa $  & $ \iterLa $ & $ \callCa $
  & $\fl$  & $\Wtime(s)$  & $\comm(s)$ 
 \\
 \hline
       4&    2268&       4
&    4313&   32344&     487
&   2.54E+14&   2.59E+04&   1.84E+02
\\
       4&    2268&       8
&    4259&   31721&     464
&   2.30E+14&   2.53E+04&   1.82E+02
\\
       4&    2268&      16
&    4218&   31230&     444
&   2.24E+14&   2.48E+04&   1.82E+02
\\
       4&    2268&      32
&    4199&   30650&     430
&   2.19E+14&   2.44E+04&   1.84E+02
\\
       4&    2268&      64
&    4199&   30506&     424
&   2.18E+14&   2.39E+04&   1.94E+02
\\
       4&    2268&     128
&    4199&   30520&     423
&   2.17E+14&   2.40E+04&   2.17E+02
\\
 \hline
 \end{tabular}
 }

%% file: Figs_fixed/costsALL_N008.tex
 {\footnotesize
 \begin{tabular}{ccc|ccc|ccc}
 \hline
  $ \Mm $  & $n_i$ & $n_0$ 
  & $ \iterNa $  & $ \iterLa $ & $ \callCa $
  & $\fl$  & $\Wtime(s)$  & $\comm(s)$ 
 \\
 \hline
       8&    1134&       8
&    4214&   37777&     482
&   8.23E+13&   1.32E+04&   2.20E+02
\\
       8&    1134&      16
&    4205&   37444&     460
&   8.07E+13&   1.30E+04&   2.22E+02
\\
       8&    1134&      32
&    4205&   37293&     438
&   7.77E+13&   1.30E+04&   2.27E+02
\\
       8&    1134&      64
&    4204&   36288&     424
&   7.35E+13&   1.29E+04&   2.35E+02
\\
       8&    1134&     128
&    4201&   35595&     422
&   7.46E+13&   1.28E+04&   2.56E+02
\\
       8&    1134&     256
&    4205&   33411&     422
&   7.22E+13&   1.26E+04&   2.89E+02
\\
 \hline
 \end{tabular}
 }

%% file: Figs_fixed/costsALL_N016.tex
 {\footnotesize
 \begin{tabular}{ccc|ccc|ccc}
 \hline
  $ \Mm $  & $n_i$ & $n_0$ 
  & $ \iterNa $  & $ \iterLa $ & $ \callCa $
  & $\fl$  & $\Wtime(s)$  & $\comm(s)$ 
 \\
 \hline
      16&     567&      16
&    4225&   34741&     457
&   2.69E+13&   6.44E+03&   2.09E+02
\\
      16&     567&      32
&    4212&   34439&     446
&   2.68E+13&   6.45E+03&   2.13E+02
\\
      16&     567&      64
&    4204&   33728&     434
&   2.63E+13&   6.52E+03&   2.21E+02
\\
      16&     567&     128
&    4202&   33550&     429
&   2.66E+13&   6.70E+03&   2.44E+02
\\
      16&     567&     256
&    4201&   32701&     423
&   2.72E+13&   7.33E+03&   2.85E+02
\\
      16&     567&     512
&    4202&   31105&     421
&   4.53E+13&   9.11E+03&   3.62E+02
\\
 \hline
 \end{tabular}
 }

%% file: Figs_fixed/costsALL_N032.tex
 {\footnotesize
 \begin{tabular}{ccc|ccc|ccc}
 \hline
  $ \Mm $  & $n_i$ & $n_0$ 
  & $ \iterNa $  & $ \iterLa $ & $ \callCa $
  & $\fl$  & $\Wtime(s)$  & $\comm(s)$ 
 \\
 \hline
      32&     283&      32
&    4222&   33819&     467
&   1.05E+13&   3.09E+03&   2.12E+02
\\
      32&     283&      64
&    4202&   33274&     454
&   1.05E+13&   3.19E+03&   2.21E+02
\\
      32&     283&     128
&    4202&   32914&     439
&   9.88E+12&   3.53E+03&   2.42E+02
\\
      32&     283&     256
&    4201&   32164&     432
&   1.69E+13&   4.46E+03&   2.84E+02
\\
      32&     283&     512
&    4203&   31110&     424
&   4.28E+13&   6.04E+03&   3.65E+02
\\
 \hline
 \end{tabular}
 }

%% file: Figs_fixed/costsALL_N064.tex
 {\footnotesize
 \begin{tabular}{ccc|ccc|ccc}
 \hline
  $ \Mm $  & $n_i$ & $n_0$ 
  & $ \iterNa $  & $ \iterLa $ & $ \callCa $
  & $\fl$  & $\Wtime(s)$  & $\comm(s)$ 
 \\
 \hline
      64&     141&      64
&    4208&   34838&     465
&   3.11E+12&   1.28E+03&   2.34E+02
\\
      64&     141&     128
&    4204&   34462&     450
&   4.96E+12&   1.61E+03&   2.57E+02
\\
      64&     141&     256
&    4205&   33665&     438
&   1.36E+13&   2.55E+03&   3.00E+02
\\
      64&     141&     512
&    4203&   32523&     431
&   4.04E+13&   4.70E+03&   3.84E+02
\\
 \hline
 \end{tabular}
 }

%% file: Figs_fixed/costsALL_N128.tex
 {\footnotesize
 \begin{tabular}{ccc|ccc|ccc}
 \hline
  $ \Mm $  & $n_i$ & $n_0$ 
  & $ \iterNa $  & $ \iterLa $ & $ \callCa $
  & $\fl$  & $\Wtime(s)$  & $\comm(s)$ 
 \\
 \hline
     128&      70&     128
&    4207&   36039&     460
&   4.97E+12&   1.28E+03&   2.72E+02
\\
     128&      70&     256
&    4205&   34975&     437
&   1.44E+13&   2.37E+03&   3.14E+02
\\
     128&      70&     512
&    4204&   33327&     435
&   4.28E+13&   4.69E+03&   3.97E+02
\\
 \hline
 \end{tabular}
 }

%% file: Figs_fixed/costsALL_N256.tex
 {\footnotesize
 \begin{tabular}{ccc|ccc|ccc}
 \hline
  $ \Mm $  & $n_i$ & $n_0$ 
  & $ \iterNa $  & $ \iterLa $ & $ \callCa $
  & $\fl$  & $\Wtime(s)$  & $\comm(s)$ 
 \\
 \hline
     256&      35&     256
&    4206&   36151&     451
&   1.53E+13&   2.18E+03&   3.28E+02
\\
     256&      35&     512
&    4204&   34242&     438
&   4.40E+13&   4.18E+03&   4.10E+02
\\
 \hline
 \end{tabular}
 }

%% file: Figs_fixed/costsALL_N512.tex
 {\footnotesize
 \begin{tabular}{ccc|ccc|ccc}
 \hline
  $ \Mm $  & $n_i$ & $n_0$ 
  & $ \iterNa $  & $ \iterLa $ & $ \callCa $
  & $\fl$  & $\Wtime(s)$  & $\comm(s)$ 
 \\
 \hline
     512&      17&     512
&    4205&   35389&     455
&   4.46E+13&   4.53E+03&   4.27E+02
\\
 \hline
 \end{tabular}
 }

%% file: Figs_adaptB/aDDM_costs.tex
 {\small 
 \begin{tabular}{c|cccc|ccc|c}
 \hline
  variant & $\Mm$ &  $\#\THm$ & $\iterNa$ & $\iterLa$ 
 & $\fl$ & $\Wtime$ & $\comm$ & $\costs$ \\
 \hline
 {\fix}, $\Mm$ =                                   
   10
 &        10 &        10 &       752 &     49122
 &     1.53E+13 &     1.99E+03 &     3.18E+01 &     2.03E+03
\\
 {\fix}, $\Mm$ =                                   
   20
 &        20 &        20 &       724 &     48455
 &     5.30E+12 &     1.53E+03 &     4.59E+01 &     1.58E+03
\\
 {\fix}, $\Mm$ =                                   
   30
 &        30 &        30 &       724 &     48276
 &     3.15E+12 &     1.13E+03 &     6.00E+01 &     1.19E+03
\\
 {\fix}, $\Mm$ =                                   
   50
 &        50 &        50 &       716 &     47810
 &     3.29E+12 &     1.04E+03 &     1.59E+02 &     1.20E+03
\\
 {\fix}, $\Mm$ =                                   
   70
 &        70 &        70 &       719 &     47951
 &     5.47E+12 &     1.31E+03 &     2.19E+02 &     1.53E+03
\\
 \hline
 {\equi}, $\kk$ =                                  
16000
 &        31 &        31 &       622 &     33903
 &     8.74E+12 &     2.11E+03 &     5.26E+01 &     2.16E+03
\\
 {\equi}, $\kk$ =                                  
12000
 &        41 &        41 &       620 &     35505
 &     6.07E+12 &     1.47E+03 &     7.16E+01 &     1.54E+03
\\
 {\equi}, $\kk$ =                                  
10000
 &        49 &        49 &       639 &     38043
 &     4.71E+12 &     1.30E+03 &     8.06E+01 &     1.38E+03
\\
 {\equi}, $\kk$ =                                  
 8000
 &        62 &        62 &       639 &     39381
 &     3.91E+12 &     1.19E+03 &     9.53E+01 &     1.28E+03
\\
 {\equi}, $\kk$ =                                  
 6000
 &        82 &        82 &       664 &     41989
 &     3.15E+12 &     1.07E+03 &     1.22E+02 &     1.19E+03
\\
 {\equi}, $\kk$ =                                  
 4000
 &       123 &       123 &       691 &     46192
 &     4.40E+12 &     1.18E+03 &     2.33E+02 &     1.41E+03
\\
 \hline
 \hline
 {\adapt}, exe                                     
    1
 &        65 &        65 &       737 &     50579
 &     2.01E+12 &     9.68E+02 &     1.24E+02 &     1.09E+03
\\
 {\adapt}, exe                                     
    2
 &        65 &        65 &       732 &     49210
 &     2.01E+12 &     9.62E+02 &     1.25E+02 &     1.09E+03
\\
 {\adapt}, exe                                     
    3
 &        53 &        53 &       734 &     49464
 &     1.95E+12 &     9.67E+02 &     1.18E+02 &     1.09E+03
\\
 {\adapt}, exe                                     
    4
 &        53 &        53 &       736 &     50241
 &     2.06E+12 &     9.81E+02 &     1.23E+02 &     1.10E+03
\\
 {\adapt}, exe                                     
    5
 &        65 &        65 &       732 &     49210
 &     2.01E+12 &     9.61E+02 &     1.25E+02 &     1.09E+03
\\
 {\adapt}, exe                                     
    6
 &        56 &        56 &       737 &     50244
 &     1.93E+12 &     9.51E+02 &     1.28E+02 &     1.08E+03
\\
 \hline
 {\adapt}, minimum 
 &        53 &        53 &       732 &     49210
 &     1.93E+12 &     9.51E+02 &     1.18E+02 &     1.08E+03
\\
 {\adapt}, average 
 &        59 &        59 &       734 &     49824
 &     1.99E+12 &     9.65E+02 &     1.24E+02 &     1.09E+03
\\
 {\adapt}, maximum 
 &        65 &        65 &       737 &     50579
 &     2.06E+12 &     9.81E+02 &     1.28E+02 &     1.10E+03
\\
 \hline
 \end{tabular}
 }

%% file: Figs_adaptK/aDDM_Kelvin_costs.tex
 {\small 
 \begin{tabular}{c|cccc|ccc|c}
 \hline
  variant & $\Mm$ &  $\#\THm$ & $\iterNa$ & $\iterLa$ 
 & $\fl$ & $\Wtime$ & $\comm$ & $\costs$ \\
 \hline
 {\fix}, $\Mm$ =                                   
   10
 &        10 &        10 &       163 &      5533
 &     7.02E+12 &     5.01E+02 &     5.46E+00 &     5.06E+02
\\
 {\fix}, $\Mm$ =                                   
   20
 &        20 &        20 &       163 &      5988
 &     2.54E+12 &     3.36E+02 &     8.00E+00 &     3.44E+02
\\
 {\fix}, $\Mm$ =                                   
   30
 &        30 &        30 &       163 &      5868
 &     1.38E+12 &     2.17E+02 &     9.10E+00 &     2.26E+02
\\
 {\fix}, $\Mm$ =                                   
   50
 &        50 &        50 &       159 &      5319
 &     6.23E+11 &     1.52E+02 &     3.73E+01 &     1.89E+02
\\
 {\fix}, $\Mm$ =                                   
   70
 &        70 &        70 &       159 &      5406
 &     5.32E+11 &     1.77E+02 &     4.50E+01 &     2.22E+02
\\
 \hline
 {\equi}, $\kk$ =                                  
16000
 &        26 &        52 &       159 &      4741
 &     3.34E+12 &     4.41E+02 &     9.17E+00 &     4.50E+02
\\
 {\equi}, $\kk$ =                                  
12000
 &        34 &        34 &       158 &      4740
 &     2.40E+12 &     3.24E+02 &     1.16E+01 &     3.36E+02
\\
 {\equi}, $\kk$ =                                  
10000
 &        41 &        41 &       159 &      5203
 &     1.92E+12 &     2.88E+02 &     1.44E+01 &     3.02E+02
\\
 {\equi}, $\kk$ =                                  
 8000
 &        51 &        51 &       160 &      5413
 &     1.38E+12 &     2.45E+02 &     1.58E+01 &     2.61E+02
\\
 {\equi}, $\kk$ =                                  
 6000
 &        68 &        68 &       161 &      5508
 &     9.64E+11 &     1.88E+02 &     2.00E+01 &     2.08E+02
\\
 {\equi}, $\kk$ =                                  
 4000
 &       102 &       102 &       162 &      5776
 &     4.47E+11 &     1.55E+02 &     2.89E+01 &     1.84E+02
\\
 {\equi}, $\kk$ =                                  
 3000
 &       136 &       136 &       162 &      5856
 &     4.01E+11 &     1.36E+02 &     3.93E+01 &     1.76E+02
\\
 {\equi}, $\kk$ =                                  
 2000
 &       140 &       140 &       159 &      5610
 &     4.33E+11 &     1.34E+02 &     4.71E+01 &     1.81E+02
\\
 \hline
 \hline
 {\adapt}, exe                                     
    1
 &        81 &        81 &       159 &      5567
 &     3.47E+11 &     1.15E+02 &     2.95E+01 &     1.45E+02
\\
 {\adapt}, exe                                     
    2
 &        83 &        83 &       159 &      5518
 &     3.37E+11 &     1.17E+02 &     2.95E+01 &     1.47E+02
\\
 {\adapt}, exe                                     
    3
 &        83 &        83 &       159 &      5501
 &     3.35E+11 &     1.14E+02 &     2.95E+01 &     1.44E+02
\\
 {\adapt}, exe                                     
    4
 &        89 &        89 &       159 &      5571
 &     3.42E+11 &     1.16E+02 &     3.11E+01 &     1.47E+02
\\
 {\adapt}, exe                                     
    5
 &        81 &        81 &       159 &      5567
 &     3.47E+11 &     1.15E+02 &     2.94E+01 &     1.45E+02
\\
 {\adapt}, exe                                     
    6
 &        87 &        87 &       159 &      5522
 &     3.36E+11 &     1.15E+02 &     2.95E+01 &     1.44E+02
\\
 \hline
 {\adapt}, minimum 
 &        81 &        81 &       159 &      5501
 &     3.35E+11 &     1.14E+02 &     2.72E+01 &     1.42E+02
\\
 {\adapt}, average 
 &        84 &        84 &       159 &      5541
 &     3.41E+11 &     1.16E+02 &     2.76E+01 &     1.43E+02
\\
 {\adapt}, maximum 
 &        89 &        89 &       159 &      5571
 &     3.47E+11 &     1.17E+02 &     2.82E+01 &     1.45E+02
\\
 \hline
 \end{tabular}
 }